\documentclass[onefignum,onetabnum]{siamart190516}
\usepackage[font={footnotesize}]{caption}
\usepackage{etoolbox}
\AtBeginEnvironment{algorithm}{\footnotesize}
\AtBeginEnvironment{algorithmic}{\footnotesize}

\usepackage [american]{babel}
\usepackage{amsmath,amsfonts,amssymb,upref}
\usepackage{units}
\usepackage{cite}
\usepackage{graphicx}        
\graphicspath{{./Figures/}}
\DeclareGraphicsExtensions {.png,.pdf,.jpg}
\usepackage{subcaption}        
\usepackage{multicol}        
\usepackage[bottom]{footmisc}
\usepackage{pgfplots}\pgfplotsset{width=0.45\textwidth,compat=1.5.1}
\usepackage{pgfplotstable}
\usepackage{booktabs}
\usepackage{array}
\usepackage{colortbl}
\usepackage[T1]{fontenc}
\usepackage[utf8x]{inputenc}
\PrerenderUnicode{äÄüÜöÖß}
\usepackage {exscale}
\usepackage{mathtools}

\usepackage{algorithm}
\usepackage{algpseudocode}
\algdef{SE}[SUBALG]{Indent}{EndIndent}{}{\algorithmicend\ }%
\algtext*{Indent}
\algtext*{EndIndent}
\algnewcommand\AlgOn[1]{\State\textbf{on} #1 \textbf{do}\Indent}
\algnewcommand\AlgEndOn{\EndIndent\State\textbf{end on}}

\newcommand{\spn}{\operatorname{span}}

\newsiamremark{remark}{Remark}
\newsiamremark{example}{Example}
\newsiamremark{hypothesis}{Hypothesis}
\crefname{hypothesis}{Hypothesis}{Hypotheses}
\newsiamthm{claim}{Claim}

\headers{Dimension-oblivious domain decomposition method}{M.~Griebel, M.~A.~Schweitzer, and L.~Troska}

\title{A dimension-oblivious domain decomposition method based on space-filling curves}

\author{Michael Griebel\thanks{Institut für Numerische Simulation, Universität Bonn, Friedrich-Hirzebruch-Allee 7, 53115 Bonn,
(\email{griebel@ins.uni-bonn.de.de}, \email{schweitzer@ins.uni-bonn.de.de}, \email{troska@ins.uni-bonn.de.de}) and Fraunhofer Institute for Algorithms and Scientific Computing (SCAI), Schloss Birlinghoven, 
53754 Sankt Augustin, Germany (\email{michael.griebel@scai.fraunhofer.de}, \email{marc.alexander.schweitzer@scai.fraunhofer.de}, \email{lukas.troska@scai.fraunhofer.de})}
\and
Marc Alexander Schweitzer\footnotemark[1]
\and
Lukas Troska\footnotemark[1]
}

\usepackage{amsopn}

\begin{document}
\maketitle

\begin{abstract} In this paper we present an algebraic dimension-oblivious two-level domain decomposition solver for discretizations of elliptic partial differential equations. The proposed parallel solver is based on a space-filling curve partitioning approach that is applicable to any discretization, i.e. it directly operates on the assembled matrix equations. Moreover, it allows for the effective use of arbitrary processor numbers independent of the dimension of the underlying partial differential equation while maintaining optimal convergence behavior. This is the core property required to attain a sparse grid based combination method with extreme scalability which can utilize exascale parallel systems efficiently. Moreover, this approach provides a basis for the development of a fault-tolerant solver for the numerical treatment of high-dimensional problems. To achieve the required data redundancy we are therefore concerned with large overlaps of our domain decomposition which we construct via space-filling curves. 
In this paper, we propose our space-filling curve based domain decomposition solver and present its convergence properties and scaling behavior. The results of numerical experiments clearly show that our approach provides optimal convergence and scaling behavior in arbitrary dimension utilizing arbitrary processor numbers.
\end{abstract}




\section{Introduction}
\label{section:introduction}
Domain decomposition (DD) techniques provide a natural way to the parallel solution of discretizations of partial differential equations (PDE). Typically, the global domain is divided into (overlapping) subdomains which are constructed geometrically. On these subdomains a linear solver is applied (in parallel) and the local independent solutions are then combined to define a global approximation. Moreover, a coarse global problem must be incorporated to achieve optimal convergence behavior. Typically, the focus in the development of DD solvers for PDEs is concerned with a minimization of the overlap of the involved subdomains to reduce the induced communication cost and to facilitate an asynchronous communication-hiding implementation. 

In this paper we however intentionally consider rather large overlaps of the involved subdomains. Our future goal is the development of fault-tolerant parallel solvers on exascale systems which require complete data redundancy to allow for data recovery when faults occur. Thus, we consider the classical domain decomposition technique in a somewhat non-standard setting of large overlap and focus on a rather algebraic construction so that it is applicable to arbitrary dimension and anisotropic meshes.
The model problem considered is a single discrete subproblem arising from the combination method discretization of a high-dimensional Poisson equation, compare Section \ref{section:motivation}. Thus, our domain decomposition approach must be dimension-oblivous and directly applicable to anisotropic meshes, i.e.~conventional DD approaches based on geometric information in a fixed dimension $d$ cannot be employed straightforwardly anymore. Instead, we utilize a space-filling curve method. Thereby, we introduce a partitioning of the high-dimensional domain without dimension-de\-pen\-dent geometric information and control the overlap of the respective subdomains via the overlap in the one-dimensional parameter space of the space-filling curve, compare Section \ref{section:domain_decomposition_method}. We consider classical iterative schemes, like the damped Richardson iteration and the correspondingly preconditioned conjugate gradient  method, to define our dimension-oblivous solver where we introduce the necessary coarse grid problem via an algebraic approach. We study their respective convergence and scaling properties. It turns out that we achieve optimal scaling and convergence behavior in arbitrary dimension with the proposed solver.


The remainder of this paper is organized as follows: In section \ref{section:motivation} we shortly give the overall motivation for our particular DD solver construction. In section \ref{section:domain_decomposition_method} we review some classical DD techniques and we introduce space-filling curves for the partitioning of grids in arbitrary dimension. Then, we present our specific construction of overlapping subdomains based on space-filling curves and introduce the resulting solver. Moreover, we give a short outlook on fault-tolerance. In section \ref{section:numerical_results} we present the results of our numerical experiments with the proposed solver in dimensions $d \leq 6$ with a focus on its convergence properties and its scaling behavior. Finally, we conclude with some remarks.

\section{Overall Motivation: The combination technique for sparse grid discretizations}
\label{section:motivation} 
First, let us shortly summarize our overall future goal which is motivated by the ever increasing number of processing cores in todays high performance computing systems. Obviously, this impressive hardware development requires a software stack that is built upon scalable parallel algorithms to allow the effective utilization of all cores. Equally important is however that resilience against hardware failures is also built into the software stack and the employed algorithms. A promising approach to the efficient and fault-tolerant utilization of supercomputers in the context of solvers for PDE especially in high dimensions is the sparse grid combination technique, see e.g. \cite{Z, BG, griebel92CombiTechnique, Griebel.Harbrecht:2014}, provided that certain mixed smoothness prerequisites are satisfied by the problem under consideration. If used in conjunction with the particular DD linear solver proposed in this paper, extreme scalability, fault-tolerance and fine-grained load balance can be attained. 

In the sparse grid discretization \cite{Z, BG} of high-dimensional elliptic PDEs via the combination method \cite{griebel92CombiTechnique, Griebel.Harbrecht:2014} 
we need to solve subproblems on grids 
with in general anisotropic mesh sizes $(h_1,\ldots h_d)$, i.e. $ h_j=2^{-l_j}, j=1,\ldots,d$, with multivariate level parameters 
$$l=(l_1,\ldots,l_d)\in \mathbb{N}^d, \mbox{ where } |l|_1:= l_1 +\ldots +l_d = L+ (d-1)-i, i=0,\ldots, d-1, l_j >0.$$
The resulting (partial) solutions $u_l(x)$, $x=(x_1,\ldots,x_d)$, are then to be combined as
\begin{equation}\label{eq:combi}
	u({ x}) \approx u^{(c)}_{L}({ x}):=
	\sum_{i=0}^{d-1} (-1)^i \left( \begin{array}{c} {d-1}\\{i}\end{array}\right) \sum_{|{l}|_1=L+(d-1)-i} u_{l}({x})
\end{equation}
to attain the final approximate solution $u^{(c)}_L$.
\begin{figure}[h]
	\begin{center}
		\includegraphics[width=0.99\textwidth]{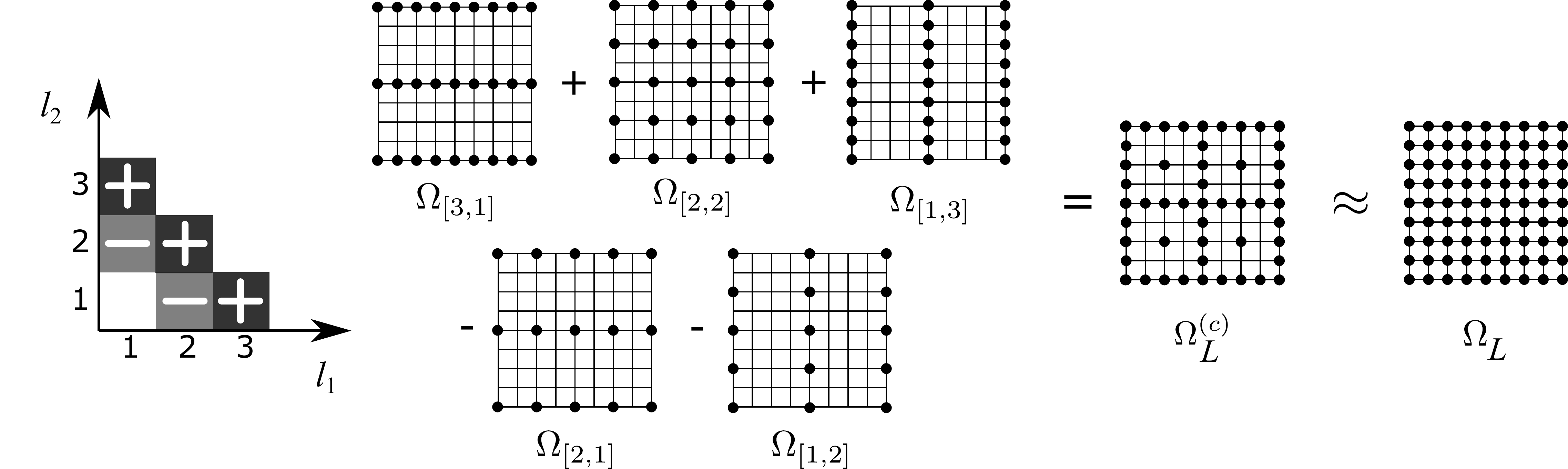}
	\end{center}
	\vspace{-1em}
	\caption{The combination method, two-dimensional case with $ L=3$.}
	\label{fig:combi_scheme}
\end{figure}
Figure~\ref{fig:combi_scheme} illustrates the construction of the combination method in two dimensions with $L=3$.
The (partial) solutions $u_l$ are completely independent of each other and thus can be computed in parallel.
Moreover, for layer $i$, we encounter 
$$ \left(  \begin{array}{c} {L+d-2-i}\\{d-1}\end{array}\right)$$ 
different subproblems, where each subproblem has approximately the same number 
\begin{equation}\label{DOFct}
N(d,i)= \prod_{j=1}^d (2^{l_j}+1)=O(2^{L+d-1-i})
\end{equation}
of degrees of freedom. We now employ an additional level of parallelism by means of the domain decomposition treatment of each of these subproblems for the partial solutions $u_l$. 
To this end, we use 
\begin{equation}\label{Pchoice}
P:= \hat P \cdot 2^{d-1-i}
\end{equation}
subdomains and thus processors/cores for each subproblem on layer $i$. This choice of a $d$- and $i$-dependent $P$ via a universal $\hat P$ in (\ref{Pchoice}) results in a subdomain size of roughly $N(d,i)/P = 2^{L+d-1-i}/P =  2^{L+d-1-i}/(\hat P \cdot 2^{d-1-i}) = 2^L/\hat P$, which is {\em independent} of $d$ and $i$. Then, in our elliptic situation and except for the coarse scale problems, only small subdomain problems of almost equal size appear. But, depending on $d$, $L$ and $\hat P$, there can be millions of these subdomains for the combination method. To be precise,  the number of subdomains for the overall set of subproblems in the combination method \eqref{eq:combi} is
$$\frac {\hat P}{(d-1)!} \sum_{k=0}^{d-1} \left( 2^{d-1-k} \prod_{i=1}^{d-1} (L+d-1-k-i)\right),$$
where each occurring subdomain problem has approximately equal size $2^L/\hat P$.

As an example, consider the goal of an overall discretization error of $10^{-12}$ in the $L^2$-norm for a problem with solution $u \in H^2_{\rm mix}(\Omega)$. Then a-priori error analysis implies that ${h\approx 10^{-6}}$ and $L=20$ are required. Moreover let us choose $\hat P=2^{10}$ in this example. For this case, Table \ref{NSPr} gives the number of subdomain problems arising in the combination method if our domain decomposition approach is employed with $P$ given by (\ref{Pchoice}) for each subproblem of layer $i$.
\begin{table}[ht]
	\caption{Overall number of subdomain problems, each one of size $2^L/\hat P = 2^{20}/2^{10}=2^{10}$, in the combination method with domain decomposition of each subproblem.} 
	\label{NSPr}
	\centering
	\begin{footnotesize}
	\begin{tabular}{c|c|c|c|c|c|c}
$d=$& 1& 2 & 3& 4& 5& 6\\
\midrule
& $1 \cdot \hat P$ & $59 \cdot  \hat P$ & $1.391 \cdot  \hat P$ & $20.889 \cdot  \hat P$ & $ 237.706 \cdot  \hat P$ & $ 1.754.744 \cdot  \hat P$\\
\end{tabular}
\end{footnotesize}
\end{table}

We see that we obtain a large number of subdomain problems to be solved in a doubly parallel way (one level of parallelism stems from the combination formula itself, the other level of parallelism stems from the domain decomposition of each subproblem of the combination method into subdomain problems). This will allow the use of extremely large parallel compute systems, i.e.~the larger $d$ is the larger the parallel system is that can be employed in a meaningful way. Furthermore the use of any fault-tolerant domain decomposition method as the inner solver for the subproblems in the combination method results in an overall fault-tolerant and parallel solver for the combination method in a natural way. There, the fault-repair mechanism is provided on the fine domain decomposition level and not just on the coarse subproblem level of the combination method, as it was previously done in \cite{Harding14, Pflueger.Bungartz.Griebel.ea.2014, Ali2016ComplexSA, Rentrop.Griebel.2020}.

Thus, an efficient parallel solver is required for each of the subproblems in the combination method. Ideally, such a linear solver needs to be applicable to anisotropic grids in arbitrary dimension while maintaining optimal convergence behavior. To this end, we propose an overlapping domain decomposition approach which is based on an (almost) algebraic construction of the overlapping subdomains via a space-filling curve approach, where we focus on larger overlaps as they are necessary to achieve data redundancy. 
{Even though this may in principle be possible to realize also with graph partitioning schemes, e.g. similar to \cite{CaiSaad1996}, which however was designed for small overlaps and small dimension, it will turn out to be trivial to attain with the proposed space-filling curve scheme indepdent of the dimension $d$.}

{Note moreover that we must simultaneously partition many of the subproblems of the combination technique with different numbers of gridpoints to utilize very large core counts and good parallel utilization altogether. 
The global load balance and scalability (for the overall combination technique solution) depends on the predictability and comparability of the computational load for the independent subproblems which are quite different due to the different anisotropic grids and consist of quite different numbers of grid points employed for each subproblem. Therefore, we will identify the number of cores $P(d,i)$ to be used for a specific subproblem with $N(d,i)$ degrees of freedom so that the load per core $\frac{N(d,i)}{P(d,i)}$ is indendent of $d$ and $i$ for each subproblem of the combination technique. Thereby we attain a good global load balance, i.e. a similar number of grid points assigned to a single core over many subproblem, while solving multiple subproblems in parallel simultaneously. This however requires that we employ a data partitioning approach that can attain a good local balance for a single subproblem of the combination technique for an arbitrary number of cores $P$ of moderate size.}

\section{Domain decomposition based on space-filling curves}
\label{section:domain_decomposition_method}
In the following we shortly review fundamental properties of some classical domain decomposition approaches which we will evaluate in our particular setting of larger overlaps before we focus on the specific construction of the involved overlapping subdomains based on space-filling curves. Then, we combine these ingredients to introduce our proposed solver and give a short outlook on its application in a fault-tolerant setting. 

\subsection{Domain decomposition techniques}
\label{subsection_section:domain_decomposition_method}
The (overlapping) domain decomposition approach is essentially a simple divide-et-impera approach to the solution of (discretized) partial differential equations and its origin can be traced back to Schwarz \cite{Schwarz}. In numerical simulations it is typically used as a preconditioner for the conjugate gradient or other Krylov iterative methods.
However, such simple domain decomposition methods can not attain fast convergence and thus, starting in the mid 80s, various techniques have been developed to introduce an additional coarse scale problem, which provides a certain amount of global transfer of information across the whole domain and thus substantially speeds up the iteration. For instance, it was shown in \cite{dryawidlund87} that the condition number of the fine grid system preconditioned by such a two-level additive Schwarz method is of the order
\begin{equation}\label{AAA}
O(1+H/\delta),
\end{equation}
where $\delta$ denotes the size of the overlap and $H$ denotes the coarse mesh size. This bound also holds for small overlap \cite{dryawidlund94} and can not be improved further \cite{brenner20}.
Thus, if the quotient of the coarse mesh size $H$ and the overlap $\delta$ stays constant, the method is indeed optimally preconditioned and weakly scalable.
For further details on domain decomposition methods see e.g. the books \cite{smithborstedgropp96, quateronivalli99, toselliwidlund04, doleanjolivetnataf15}.

We obtain a two-level additive Schwarz method as follows: Consider an elliptic differential equation 
\begin{equation}\label{ellprob}
\mathcal{L} u=f \text{ in } \Omega \subset \mathbb{R}^d,
\end{equation}
e.g. the simple Poisson problem on a $d$-dimensional cube. Using a conforming finite element method, a finite difference approach or a finite volume discretization involving $N$ degrees of freedom and mesh size $h\approx N^{-1/d}$ (for the ease of notation we consider an isotropic mesh for the discussion of the fundamental properties of the domain decomposition technique), we arrive at the system of linear equations
\begin{equation}\label{LES}
	A x = b
\end{equation}
with sparse  stiffness matrix $A \in \mathbb{R}^{N \times N}$, right hand side vector $b \in \mathbb{R}^N$ and unknown coefficient vector $x \in \mathbb{R}^N$, which needs to be solved.
Now suppose that 
$$\Omega=\bigcup_{i=1}^P \Omega_i$$ 
is covered by a finite number $P$ of well-shaped subdomains $\Omega_i$ of diameter $\approx H$ which might locally overlap. It is silently assumed that $h\ll H$ and that the subdomains are aligned with the fine mesh. 
Moreover denote by $N_i$ the number of grid points associated to each $\Omega_i$, i.e.~the degrees of freedom associated to the subdomains $\Omega_i, i=1,\ldots,P$.
Then denote by $R_i: \mathbb{R}^{N} \to \mathbb{R}^{N_i}$ the restriction operators, which map the entries of the coefficient vector $x$ corresponding to the full grid on $\Omega$ to the coefficient vectors $x_i$ corresponding to the local grids on the subdomains $\Omega_i$.
Analogously denote by $R^T_i: \mathbb{R}^{N_i} \to \mathbb{R}^{N}$ the extension operators, which map the coefficient vectors from the local grid on the subdomains $\Omega_i$ to that of the full grid on $\Omega$ via the natural extension by zero.
Then the local stiffness matrices associated to the subdomains $\Omega_i$ can be denoted as $A_i \in \mathbb{R}^{N_i \times N_i}$  with ${A_i := R_i A R^T_i}$ and the resulting one-level additive Schwarz operator is defined as
\begin{equation}\label{AS2a}
C_{(1)}^{-1}:=   \sum_{i=1}^P  R_i^T A_i^{-1} R_i.
\end{equation}
Finally, we add an additional coarse space problem with dimension $N_0$ as a second level via the restriction operator $R_0: \mathbb{R}^{N} \to \mathbb{R}^{N_0}$, which maps from the full grid on $\Omega$ to a respective global coarse grid. Similarly as for the subdomain problems, the associated coarse stiffness matrix $A_0$ can be generated via the Galerkin approach as $A_0 := R_0 A R^T_0$. Altogether, with the one-level additive Schwarz operator \cref{AS2a}
we obtain the two-level additive Schwarz operator
\begin{equation}\label{AS2}
	C_{(2)}^{-1}:= R_0^T A_0^{-1} R_0 + C_{(1)}^{-1} = \sum_{i=0}^P  R_i^T A_i^{-1} R_i.
\end{equation}
These additive Schwarz operators \eqref{AS2a} and \eqref{AS2} can then be used as preconditioners in Krylov methods, like the conjugate gradient method, or with the help of an appropriate scaling parameter, they can be used directly to define Richardson-like linear iterative solvers.
A notational variant based on space splittings is given in \cite{Griebel.Oswald:2019}.
Note here that there are more sophisticated space splittings which follow the Bank-Holst technique \cite{BH2003},
where the coarse problem is formally avoided by including a redundant copy of it into each of the subdomain problems with $i=1,\ldots,P$. We will indeed follow this approach later on.

Now, if the condition number $\kappa(C_{(2)}^{-1}A)= \lambda_{\max}(C_{(2)}^{-1}A)/\lambda_{\min}(C_{(2)}^{-1}A)$ of the preconditioned system is
independent of the number $P$ of subproblems for fixed $N$, we obtain
{\em strong scalability}. If it is independent of $P$ while the quotient $N/P$ is fixed, i.e.~the problem size per subdomain and thus per processor stays fixed, we obtain {\em weak scalability}. 
Moreover, if it is independent of the number $N$ of degrees of freedom, we would have an optimally preconditioned method, which however still may depend on $P$ and might thus not be scalable.
Note furthermore that we employ here for reasons of simplicity a direct solver for $A_i^{-1}$ on all subdomains and for $A_0^{-1}$ on the coarse scale, which involves Gaussian elimination and comes with a certain cost. However, the corresponding matrix factorization needs to be performed just once at the beginning and, in the plain linear iteration or in the preconditioned conjugate gradient iteration, only the cheaper backward and forward steps need to be employed. Alternatively, approximate iterative methods might be used as well, like the multigrid or BPX-multilevel method, which would result in optimal linear cost for the subproblem solves. This given, to achieve a mesh-independent condition number for the preconditioned system $C_{(2)}^{-1}A$ with $C_{(2)}$ (\ref{AS2}), one usually chooses for the coarse problem a suitable finite element space on the mesh of domain partitions, where a linear FE space will do for a second-order elliptic problem such as (\ref{ellprob}). Mild shape regularity assumptions on the overlapping subdomains $\Omega_i$ then guarantee robust condition number estimates of the form  
$
\kappa\le c(1+\frac H \delta),
$
see \cite[Theorem 3.13]{toselliwidlund04}.
Dropping the coarse grid problem, i.e.~considering a one-level preconditioner as in (\ref{AS2a}) without the coarse problem $ R_0^T A_0^{-1} R_0 $,
yields the worse bound 
$\kappa\le c H^{-2}(1+H/\delta)$.
Note that, even though these estimates imply a deterioration of the condition number proportional to $\delta^{-1}$ if $\delta\to h$, in practice good performance has been observed when only a few layers of fine mesh cells form the overlap near the boundaries of the $\Omega_i$. 
With the use of an additional coarse grid problem based on piecewise linears, an optimal order of the convergence rate is then guaranteed for elliptic problems. The additional coarse grid problem results in a certain communication bottleneck,  which however can by principle not be avoided and is inherent in all elliptic problems.  Similar issues arise for multigrid and multilevel algorithms as well, but these methods are more complicated to parallelize in an efficient way on large compute systems. Moreover their achieved convergence rate and cost complexity is not better than for the domain decomposition approach with coarse grid, at least in order, albeit the order constant might be.


Apart from a mesh-based geometric coarse grid problem, we can also derive a coarse problem in a purely {\em algebraic} way. To this end,  let $V_H$ be a coarse space with $N_0:=\dim(V_H)$ and let $Z$ be a basis for it, i.e.~$V_H= \spn Z$. The restriction $R_0:\mathbb{R}^{N}\to \mathbb{R}^{N_0}$ from the fine grid space $V_h$ to the coarse space $V_H$ can be (algebraically) defined as the matrix $R_0:=Z^T\in \mathbb{R}^{N_0 \times N}$ and the coarse space discretization is again given via the Galerkin approach as
$A_0 := R_0 A R^T_0.$

There are now different specific choices for $Z$ or $R_0^T$, respectively.
In \cite{nicolaides87} it was suggested to employ the kernel of the underlying differential operator as coarse space, i.e.~the constant space for the (weak) Laplace operator. Thus, with $N_0=P$,
$$R_0^T:=(R_i^T D_iR_i {\bf 1})_{1 \leq i \leq P}$$
with $ {\bf 1}=(1,\ldots,1)^T$ and with diagonal matrices $D_i$ such that a partition of unity results, i.e.
$$\sum_{i=1}^{P} R_i^T D_iR_i = I.
$$ 
Indeed, it is observed that the associated two-level preconditioner gives good results and weak scaling is achieved in practice. Note that this approach can be easily generalized to the case $N_0 = q P$ such that $q$ degrees of freedom are associated on the coarse scale to each subdomain instead of just one.
Improved variants, namely the balanced method \cite{mandel93} and the deflation method \cite{vuiknabben06}, had been developed subsequently.
With the definitions
\begin{equation}\label{balan}
 F:=R_0^T A_0^{-1}R_0 \quad \mbox{ and } \quad G:=I-AF
\end{equation}
and the plain one-level additive Schwarz operator $C_{(1)}^{-1}$ from (\ref{AS2a}),
we get the additively corrected operator $C_{(1)}^{-1}+F$ due to Nicolaides \cite{nicolaides87}, the deflated approach $G^TC_{(1)}^{-1} +F$ \cite{vuiknabben06} and the balanced version \cite{mandel93}
\begin{equation}\label{bala}
	C_{(2),bal}^{-1}:= G^TC_{(1)}^{-1}G+F.
\end{equation}
Closely related are agglomeration techniques inspired by the alge\-braic multi\-grid meth\-od, volume agglomeration methods stemming from multigrid in the context of finite volume discretizations and so-called spectral coarse space constructions, see e.g. \cite{doleanjolivetnataf15}.

\subsection{Space-filling curves}
\label{subsection:space_filling_curves}
A main question for domain decomposition methods is how to construct the partition $\{ \Omega_i\}_{i=1}^P$ in the first place. To this end, for a fixed number $P$, one aim is surely to derive a set of subdomains which involves equal size and thus equal computational load for each of the subdomains. If we just consider uniform fine grid discretizations, a simple uniform and geometric decomposition of the mesh would need $P=\bar P^d$ with $\bar P$ being the amount of subdomain splittings in each coordinate direction. This however prohibits the use of a slowly growing number of processors $P$ especially in higher dimensions $d$.
{Similar issues arise for anisotropic mesh sizes $h=(h_1,\ldots, h_d)$ with $h_i \neq h_j$ as they originate in the combination method, more general quasi-uniform finite element meshes or adaptive meshes when using geometric partitioning schemes based on coordinate bisection. One traditional approach to overcome this issue is the use of graph partitioning schemes which however have not been studied in high dimensions so far. Thus, we follow the lead of \cite{griebelzumbusch00,griebelzumbusch01} in this paper and rely on {\em space-filling curves} to construct the partition $\{ \Omega_i\}_{i=1}^P$. 
This way, the overall number $N$ of degrees of freedom is partitioned for $P$ processors into $N/P$-sized subproblems regardless of the dimension $d$ of the PDE and the number of available processors $P$. For further details see also \cite{GriebelKnapekZumbusch2007, Bader2013, zumbusch03}.}

An important aspect of any space-filling curve is its locality property, i.e.~points that are close to each other in $[0,1]$ should tend to be close after mapping into the higher dimensional space. But how close or apart will they be?
To this end, let us consider the Hilbert curve. In \cite{zumbusch03} it was shown that the $d$-dimensional Hilbert curve $s(x)$ is H\"older continuous, but nowhere differentiable and that we have, for any $x,y \in [0,1]$, the H\"older condition
	\[	 |s(x)-s(y)|_2 \leq c_s |x-y|^{1/d} \quad \quad \mbox{ with } c_s = 2 \sqrt{d+3}. 	\]
The exponent $1/d$ of the H\"older condition is optimal for space-filling curves and cannot be improved. 

Originally, the Hilbert curve and space-filling curves in general were defined in two dimensions only. Their recursive construction, however, can be generalized to higher space dimensions $d$, i.e.~to curves  $s:[0,1]\to[0,1]^d$.
There exist codes for a number of space-filling curves, and especially for the Hilbert curve e.g. on \texttt{github} \cite{chernoch}. A recursive pseudo-code for $d$ dimensions is e.g. given in \cite{zumbusch03}. 
The approach in \cite{Butz1971} was further developed 
in e.g.~\cite{Moore2000, Skilling2004}, we use \cite{adishavit} as the basis for our implementation. 

\subsection{Proposed preconditioner}
\label{subsection:proposed_perconditioner}
Now we discuss the main features of our domain decomposition method for $d$-dimensional elliptic problems $\mathcal{L} u=f$.
For reasons of simplicity, we consider here the unit domain $\Omega=[0,1]^d$ and employ Dirichlet boundary conditions on $\delta \Omega$. 
The discretization is done with a uniform but in general anisotropic mesh size $h=(h_1,\ldots,h_d)$, where $h_j=2^{-l_j}$ with multivariate level parameter $l=(l_1,\ldots,l_d)$, which gives the global mesh $\Omega_h$. The number of interior grid points, and thus  the number of degrees of freedom, is then 
\begin{equation}\label{dof}
	N:=\prod_{j=1}^d (2^{l_j}-1).
\end{equation}

Now, for the discretization of $\mathcal{L}u=f$ on the grid $\Omega_h$, we employ the simple finite difference method (or the usual finite element method with piecewise $d$-linear basis functions) on $\Omega_h$ which results in the system of linear equations $A x = b$ with sparse system matrix $A \in \mathbb{R}^{N \times N}$ and right hand side vector $b \in \mathbb{R}^N$. 

\begin{figure}[tb]
\centering
  \includegraphics[width=0.3\textwidth,height=0.3\textwidth]{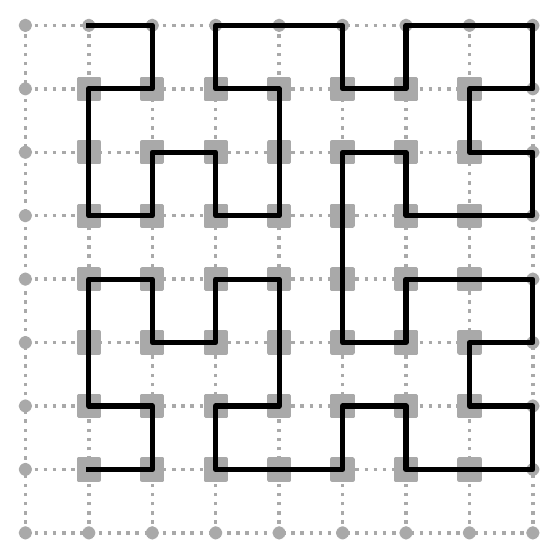}\hfill
  \includegraphics[width=0.3\textwidth,height=0.3\textwidth]{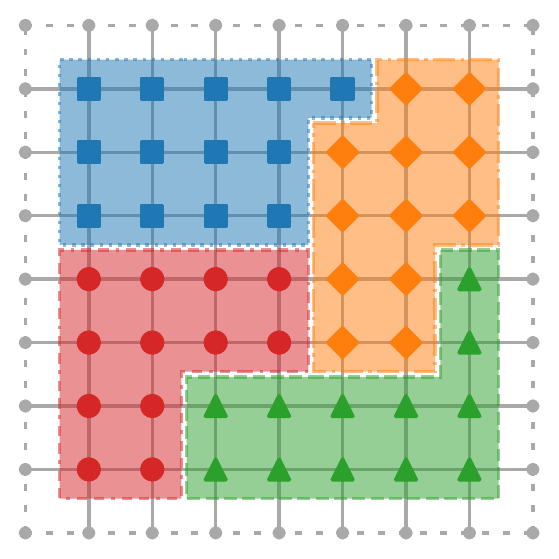}\hfill
  \includegraphics[width=0.3\textwidth,height=0.3\textwidth]{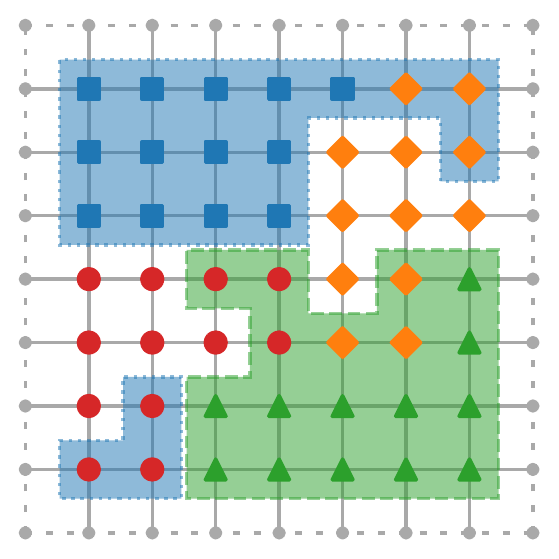}
	\caption{Decomposition of an isotropic grid with ${l}=(3,3)$ by the Hilbert curve approach (left) and corresponding disjoint subdomains (middle). Construction of two overlapping subdomains (square/blue and triangle/green) by enlargement using the Hilbert space-filling curve (right).}
\label{sfcmap}
\label{sfccmap_overlap}
\end{figure}

\begin{figure}[tb]
\centering
  \includegraphics[width=0.3\textwidth,height=0.3\textwidth]{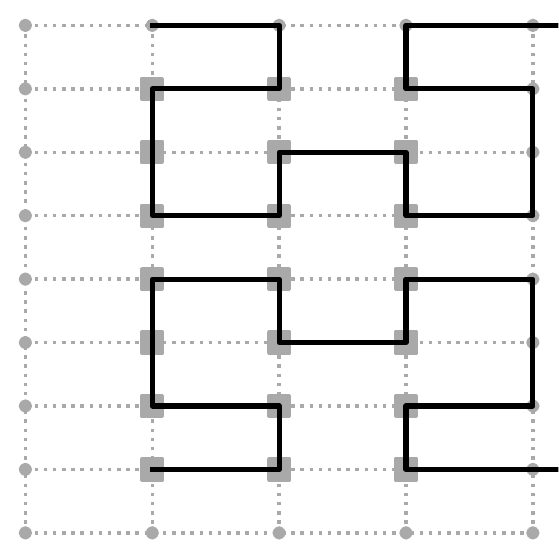}\hfill
  \includegraphics[width=0.3\textwidth,height=0.3\textwidth]{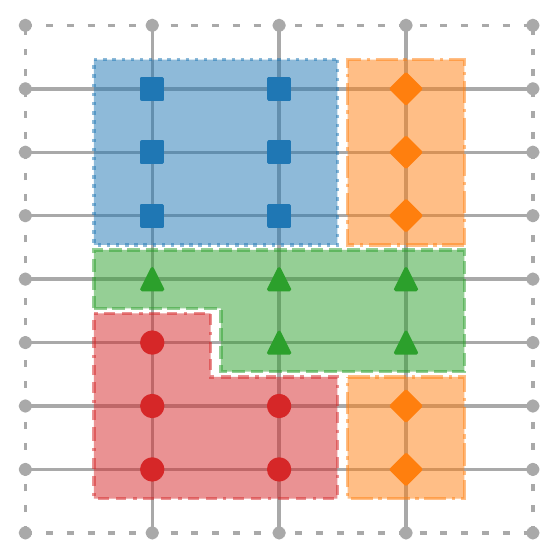}\hfill
  \includegraphics[width=0.3\textwidth,height=0.3\textwidth]{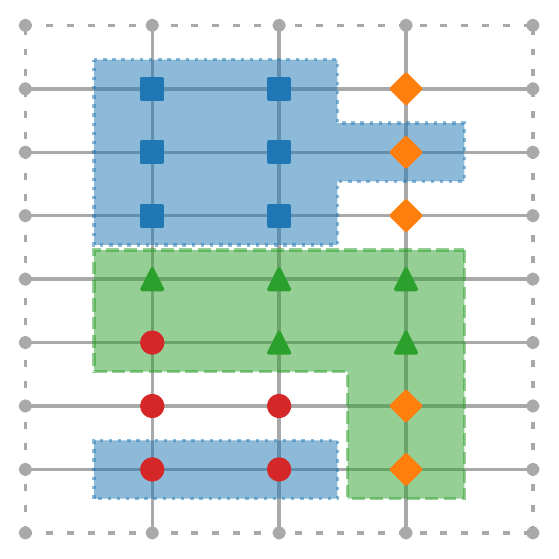}
\caption{Decomposition of an anisotropic grid with ${l}=(2,3)$ by the Hilbert curve approach (left) and corresponding disjoint subdomains (middle). Construction of two overlapping subdomains (square/blue and triangle/green) by enlargement using the Hilbert space-filling curve (right).}
\label{sfcmap2}	
\label{sfccmap_overlap2}
\end{figure}

Next, we consider the case of $P$ subdomains.  To generate an initial partition of $P$ disjoint subdomains $\{\Omega_i\}_{i=1}^P$ of equal size, we employ the space-filling curve approach and in principle just map our $d$-dimensional (interior) grid points $x_k \in \Omega_h$ by means of the inverse discrete space-filling curve $s_n^{-1}$ with sufficiently large $n$ to the one-dimensional unit interval $[0,1]$. Then we simply partition the one-dimensional, totally ordered sequence of $N$ points into a consecutive one-dimensional set of disjoint subsets of approximate size $ N/P$ each. To this end, we first determine the remainder $r:=N- P \lfloor N/P\rfloor$. This gives us the number $r$ of subdomains which have to possess $\lfloor N/P\rfloor +1$ grid points, whereas the remaining $P-r$ subdomains possess just $\lfloor N/P\rfloor $ grid points. Thus, with 
\begin{equation}\label{sloc}
	\tilde N_i:=\lfloor N/P\rfloor +1, i=1,\ldots, r \mbox{ and } \tilde N_i=\lfloor N/P\rfloor, i=r+1, \ldots,P,
\end{equation}
we assign the first $\tilde N_1$ points to the set $\tilde \Omega_1$, the second $\tilde N_2$ points to the set $\tilde \Omega_2$, and so on. Since the $\tilde N_i$ differ at most by one, we obtain a perfectly balanced partition $\{\tilde \Omega_i\}_{i=1}^P$. The basic partitioning approach by means of the Hilbert curve is shown for the two-dimensional case in Figures \ref{sfcmap} and \ref{sfcmap2}. Note here that the resolution of the discrete isotropic space-filling curve is chosen as the one which belongs to the largest value $\max_{j=1,\ldots,d} l_j$ of the entries of the level parameter $l=(l_1,\ldots,l_d)$, i.e.~to the finest resolution in case of an anisotropic grid. 


In a next step, we enlarge the initial disjoint subdomains $\tilde \Omega_i$, i.e.~the corresponding subsets of grid point indices, in a specific way to create overlap.
This is not done as in conventional, geometry-based overlapping domain decomposition methods by adding a $d$-dimensional mesh-stripe with diameter $\delta$ of grid cells at the boundary of the $d$-dimensional subdomains that are associated to the sets $\tilde \Omega_i$. Instead, we deliberately stick to the one-dimensional setting which is induced by the space-filling curve: We choose an overlap parameter ${\gamma \in \mathbb{R}, \gamma >0,}$ and enlarge the index sets $\tilde \Omega_i$ as
\begin{equation}\label{enlarge}
	\Omega_i := \tilde \Omega_i \cup \bigcup_{k=1}^{\lfloor \gamma \rfloor} \left(\tilde \Omega_{i-k} \cup \tilde \Omega_{i+k} \right) \cup \tilde \Omega_{i-\lfloor \gamma \rfloor-1}^{\eta,+} \cup \tilde \Omega_{i+\lfloor \gamma \rfloor+1}^{\eta,-}.
\end{equation}
Here, with $\eta:=\gamma-\lfloor \gamma\rfloor$, the set $\tilde \Omega_{k}^{\eta,+}$ is the subset of $\tilde \Omega_{k}$ which contains its last $\lceil \eta N_{k}\rceil$ indices with respect to the space-filling curve ordering, while the set $\tilde \Omega_{k}^{\eta,-}$ is the subset of $\tilde \Omega_{k}$ which contains its first $\lfloor \eta N_{k}\rfloor$ indices.
For example, for $\gamma=1$, we add to $\tilde \Omega_i$ exactly the two neighboring index sets $\tilde \Omega_{i+1}$ and $\tilde \Omega_{i-1}$, for $\gamma=2$,  we add the four sets $\tilde \Omega_{i+1},\tilde \Omega_{i+2} $ and  $\tilde \Omega_{i-1},\tilde\Omega_{i-2}  $. 
For $\gamma=0.5$ we would add those halves of the indices of $\tilde \Omega_{i+1}$ and $\tilde \Omega_{i-1}$ which are the ones closer to $\tilde \Omega_i$, etc. 
Moreover, to avoid any special treatment for the first and last few  $\tilde \Omega_i, i=1,2,..  $ and $i=P,P-1,..$, we cyclically close the enumeration of the subsets, i.e.~the left neighbor of $\tilde \Omega_1$ is set as $\tilde \Omega_P$ and the right neighbor of $\tilde \Omega_P$ is set as $\tilde \Omega_1$. Note that, besides $\gamma$, also the specific values of the $\tilde N_i$ enter here. Examples of the enlargement of the index sets from $\tilde \Omega_i$ to $\Omega_i$ with $\gamma=0.25$ are given in Figures \ref{sfccmap_overlap} and \ref{sfccmap_overlap2}. Here we show the induced domains $\Omega_1$ and $\Omega_3$ only.

This way, an overlapping partition $\{\Omega_i\}_{i=1}^P$ is finally constructed. Note at this point that, depending on $N,P, \gamma$ and the respective space-filling curve type, each subdomain of the associated decomposition is not necessarily fully connected, i.e.~there can exist some subdomains with point distributions that geometrically indicate some further separation of the subdomain, see e.g. $\Omega_4 $ (red) in Figure \ref{sfccmap_overlap}. But since we merely deal with index sets and not with geometric subdomains, this causes no practical issues. 
Note furthermore that there is not necessarily always a complete overlap but sometimes just a partial overlap between two adjacent subdomains being created by our space-filling curve approach.
But this also causes no practical issues. Moreover, in contrast to many other domain decomposition methods, where a goal is to make the overlap as small as possible, our approach rather aims at larger overlaps along the space-filling curve, which will yield the sufficient redundancy of stored data that is needed to attain fault tolerance.
 
Finally, we set up our coarse space problem. To this end, the size $N_0$  of the coarse problem is given via the number $P$ of subdomains and the number of degrees of freedom $n_i$ considered on the coarse scale per subdomain, i.e.
$N_0 := \sum_{i=1}^P n_i.$
If we let all $n_i$ be equal to a fixed integer $q\in \{1,\ldots, \lfloor \frac N P \rfloor\}$, i.e.~$n_i=q$, then $N_0= q \cdot P$. The mapping from the fine grid to the coarse space is now given by means of the restriction matrix $R_0 \in  \mathbb{R}^{N_0 \times N}$ and its entries. Again, we avoid any geometric information here, i.e.~we do not follow the conventional approach of a coarse grid with associated coarse piecewise linear basis functions. Instead we derive the coarse scale problem in a purely {\em algebraic} way.
For that, we resort to the non-overlapping partition $\{\tilde \Omega_i\}_{i=1}^P$ and assign the values of zero and one to the entries of the restriction matrix as follows:
Let $q$ be the constant number of coarse level degrees of freedom per subdomain. With $\tilde N_i=\lfloor N/P\rfloor+1$ if $i \leq (N \mod P) $ and $\tilde N_i=\lfloor N/P\rfloor$ otherwise, which denotes the size of $\tilde \Omega_i$, we have the index sets $\tilde \Omega_i =\left\{ \sum_{j = 1}^{i-1} \tilde N_j + 1, \ldots, \sum_{j=1}^i \tilde N_j\right\} $. Now let $\tilde \Omega_{i, m}$ be the $m$-th subset of $\tilde \Omega_i$ with respect to the size $q$ in the space-filling curve ordering, 
i.e.
$$\tilde \Omega_{i, m} = \left\{\sum_{j=1}^{i-1} \tilde N_j + \sum_{n=1}^{m-1} \tilde N_{i, n} + 1, ..., \sum_{j=1}^{i-1} \tilde N_j + \sum_{n=1}^m \tilde N_{i, n}\right\},
$$
with associated size $\tilde N_{i,m}$ for which, with $q$ coarse points per domain, we have
\begin{equation*}
	\tilde N_{i, m} = \left( 
	\begin{array}{ll}
		\lfloor \tilde N_i/q \rfloor + 1  & \mbox{if } m \leq (\tilde N_i \mod q), \\ 
		\lfloor \tilde N_i/q \rfloor 	 & \mbox{otherwise}.
	\end{array}
\right.
\end{equation*}
Then,  for $i=1,...,P, m=1,...,q, j=1,...,N$, we have
\begin{equation}\label{restr}
	(R_{0})_{((i-1)q+m,j)}= \left( 
	\begin{array}{ll}
		1 & \mbox{if } j \in \tilde \Omega_{i,m}, \\ 	0 & \mbox{otherwise}.
	\end{array}
\right.
\end{equation}
In this way, a basis is implicitly generated on the coarse scale: Each basis function is constant in the part of each subdomain of the non-overlapping partition which belongs to the
$\tilde \Omega_{i,m}$, where $q$ piecewise constant basis functions with support on $\tilde \Omega_{i,m}$ are associated to each $\tilde \Omega_i$.

The coarse scale problem is then set up via the Galerkin approach as 
\begin{equation}\label{csm}
A_0 = R_0 A R^T_0.
\end{equation}
Here we follow the Bank-Holst technique \cite{BH2003} and store a redundant copy of it on each of the $P$ processors together with the respective subdomain problem.  This way, the coarse problem is formally avoided. Moreover the coarse problem is redundantly solved on each processor. It interacts with the respective subproblem solver in an additive, parallel way, i.e.~we solve the global coarse scale problem and the local fine subdomain problem independently of each other on the processor.

Finally, we have to deal with the overcounting of unknowns due to the overlapping of fine grid subdomains. To this end, we resort to the partition of unity on the fine level
\begin{equation}\label{pum} 
	I = \sum_{i=1}^P R_i^T D_i R_i 
\end{equation}
with properly chosen diagonal matrices $D_i \in \mathbb{R}^{N_i \times N_i}$ to define the weighted one-level operator 
\begin{equation}\label{w1levprec}
C_{(1),D}^{-1} := \sum_{i=1}^P  R_i^T D_i A_i^{-1} R_i.
\end{equation}
and the respective two-level domain decomposition operator
\begin{equation}\label{addSW}
	C_{(2),D}^{-1} :=  R_0^T A_0^{-1} R_0 + C_{(1),D}^{-1} = \sum_{i=0}^P  R_i^T D_i A_i^{-1} R_i
\end{equation}
where $D_0:=I$. In complete analogy to \eqref{bala} we then introduce the balanced version

\begin{equation}\label{balSW}
C_{(2),D,bal}^{-1}:= G^T C_{(1),D}^{-1} G + F.
\end{equation}

There are different choices of the $D_i$'s since the above condition (\ref{pum}) does not have a unique solution.
A natural choice is obviously
\begin{equation}\label{Di}
	(D_i)_{{j,j}} = 1/ |\left\{\Omega_{i'}, i'=1,\ldots P : j \in \Omega_{i'}\right\} |,
\end{equation}
which locally takes for each index $j\in \Omega_i$ the number of domains that overlap this index into account. Note that, for such general diagonal matrices $D_i$,  the associated 
$C_{(2),D}$ is not a symmetric operator. 
If however each $D_i$ is chosen as $\omega_i I_i$ with a scalar positive weight $\omega_i$, then, on the one hand, the partition of unity property is lost, but, on the other hand, symmetry is regained and we have a weighted, two-level overlapping additive Schwarz operator, which can also be used as a preconditioner for the conjugate gradient method. To this end, a sensible choice is
\begin{equation}\label{omegai}
	\omega_i =\max_j (D_i)_{(j,j)}
\end{equation}
with $D_i$ from (\ref{Di}). Note furthermore that due to our particular construction of the overlapping partition via space-filling curves there is a much more natural choice of $D_i$, in particular $\omega_i$, for certain choices of the overlap parameter $\gamma$, see section \ref{subsection:outlook_fault_tolerance} below.

In any case, we obtain the two-level Schwarz operator \eqref{addSW}, with its application given by Algorithm \ref{Alg:DDMAlgorithm2}, where the setup phase, which is executed just once before the iteration starts, is described in Algorithm \ref{Alg:DDMAlgorithm}. The application of $C_{(2),D,bal}^{-1}$ can be described in a similar fashion. The convergence properties of the associated preconditioned damped linear/Richardson-type iteration (Algorithm \ref{Alg:Richardson}) are well known. In fact it is a special case of the additive Schwarz iteration as studied in \cite{Griebel.Oswald:1995*1}. There, in Theorem 2, it is stated that an additive method, e.g. induced by (\ref{addSW}) and damped by a scalar $\xi$, indeed converges in the energy norm for any damping $0 < \xi  <2/\lambda_{\max}$ with the convergence rate 
$$\rho =\max\{|1-\xi \lambda_{\min}|,|1-\xi \lambda_{\max}|\},
$$
where $\lambda_{\min}$ and $\lambda_{\max}$ are the smallest and largest eigenvalues of $C_{(2),D}^{-1}A$, provided that $C_{(2),D}^{-1}$ is a symmetric operator, which for example is the case for the choice (\ref{omegai}). Moreover, the optimal damping $\xi^*$ and optimal convergence rate $\rho^*$ are given by 
\begin{equation}\label{redrate}
\rho^* = 1- \frac 2 {1+\kappa}  \text{ with } \kappa =\frac{\lambda_{\max}}{\lambda_{\min}}, \quad \xi^* = \frac 2 {\lambda_{\min}+\lambda_{\max}}.
\end{equation} 
The proof is exactly the same as for the conventional Richardson iteration. 
To this end, the numbers $\lambda_{\max}$ and $\lambda_{\min}$ need to be explicitly known to have the optimal damping parameter $\xi^*$, which is of course not the case in most practical applications. Then good estimates, especially for $\lambda_{\max}$, are needed to still obtain a convergent iteration scheme, albeit with a somewhat reduced convergence rate. Note at this point that for the general non-symmetric case, i.e.~for the general choice (\ref{Di}), this convergence theory does not apply. In practice, however, convergence can still be observed.

\begin{algorithm}[hbt!]
\caption{Overlapping two-level additive Schwarz iteration with space-filling curve: Setup phase.}
\label{Alg:DDMAlgorithm}

\begin{algorithmic}[1]
\AlgOn{every processor $i=1,\ldots,P$}
    \State \parbox[t]{\dimexpr\linewidth-\algorithmicindent-\algorithmicindent}{Set input parameters: $d$, $l=(l_1,\ldots,l_d)$, $P$, $\gamma$, $q$, type of space-filling curve.\strut}
    \State \parbox[t]{\dimexpr\linewidth-\algorithmicindent-\algorithmicindent}{Derive $N$ from (\ref{dof}), set $\tilde N_i, i=1,\ldots,P$ as in (\ref{sloc}) and set $N_0=q\cdot P$.\strut}
    \State \parbox[t]{\dimexpr\linewidth-\algorithmicindent-\algorithmicindent}{Compute the index vector $sfc\_index$ of length $N$ from the $d$-dimensional grid point indices $k=(k_1,\ldots,k_d)$, $k_j=1,\ldots,2^{l_j}-1$, $j=1,\ldots d,$ according to the space-filling curve by means of {\em cmp}$((k_1, ..., k_d), (k'_1, ..., k'_d))$.\strut}
    \State \parbox[t]{\dimexpr\linewidth-\algorithmicindent-\algorithmicindent}{Derive the disjoint subdomain index sets $\{\tilde \Omega_i\}_{i=1}^P$ by splitting the overall index set into $P$ subsets $\tilde\Omega_i$ of consecutive indices, each of size $\tilde N_i$. This is simply done by storing two integers $\tilde ta_i,\tilde tb_i$, which indicate where the local index sequence of $\tilde \Omega_i$ starts and ends in $sfc\_index$.\strut}
    \State \parbox[t]{\dimexpr\linewidth-\algorithmicindent-\algorithmicindent}{Derive the overlapping subdomain index sets $\{\Omega_i\}_{i=1}^P$ by enlarging the $\tilde\Omega_i$ with $\gamma$ as in (\ref{enlarge}). Again, this is simply done by storing two integers  $ ta_i,tb_i$, which indicate where the local index sequence of $\Omega_i$ starts and ends in $sfc\_index$.\strut}
    \State \parbox[t]{\dimexpr\linewidth-\algorithmicindent-\algorithmicindent}{Set up a map to neighboring grid points that are not in $\Omega_i$, i.e.~store their global indices, to later determine the column entries of the stiffness matrix that are situated outside of $\Omega_i$.\strut}
    \State \parbox[t]{\dimexpr\linewidth-\algorithmicindent-\algorithmicindent}{Set the rows of $A$ that belong to $\Omega_i$, i.e.~store the rows of $A$ with indices $j \in [ta_i,tb_i]$ in CRS format.\strut}
    \State \parbox[t]{\dimexpr\linewidth-\algorithmicindent-\algorithmicindent}{Initialize the part of the starting iterate $x^0$ and the part of $b$ that belong to $\Omega_i$.\strut}
	\State \parbox[t]{\dimexpr\linewidth-\algorithmicindent-\algorithmicindent}{Derive the rows of the matrix $R_0$ from (\ref{restr}) with indices $j \in [(i-1)q+1, \dots,iq] $ and store them in CRS format.\strut}
	\State \parbox[t]{\dimexpr\linewidth-\algorithmicindent-\algorithmicindent}{Compute the rows of the coarse scale matrix $A_0$ as in (\ref{csm}) that belong to $\Omega_i$,i.e.~with the indices $j \in [(i-1)q+1, \dots,iq] $, and store them in CRS format.\strut}
\AlgEndOn
\end{algorithmic}
\end{algorithm}

\begin{algorithm}[hbt!]
\caption{Overlapping two-level additive Schwarz iteration with space-filling curve: Application of operator.}
\label{Alg:DDMAlgorithm2} 
\begin{footnotesize}
\hspace*{\algorithmicindent} \textbf{Input:} Vector $g$ \\
\hspace*{\algorithmicindent} \textbf{Output:} Vector $h\coloneqq C_{(2),D}^{-1} g$ 
\begin{algorithmic}[1]
\AlgOn{every processor $i=1,\ldots,P$}
\State \parbox[t]{\dimexpr\linewidth-\algorithmicindent}{Compute the part $g_i=R_i g$ of the vector $g$ that belongs to $\Omega_i$.\strut}
\State Solve the local subproblems  $A_i d_i=g_i$
\State Solve redundantly the coarse scale problem $A_0 d_0 = R_0 r$.
\State \parbox[t]{\dimexpr\linewidth-\algorithmicindent}{Compute the part $h_i=R_ih$ of the vector $h= \sum_{i=0}^P R_i^T D_i d_i$ that belongs to $\Omega_i$.\strut}
\AlgEndOn
\end{algorithmic}
\end{footnotesize}
\end{algorithm}

\begin{algorithm}[hbt!]
\caption{Preconditioned damped linear/Richardson iteration.}
\label{Alg:Richardson} 

\begin{algorithmic}[1]
\State Set k=1.
\While{not converged}
\State Compute the residual $r^k=b-Ax^k$
\State Compute $h^k$ by applying a preconditioner to $r^k$ (e.g.\ \eqref{addSW} via Algorithm \ref{Alg:DDMAlgorithm2} or \eqref{balSW}).
\State Update the iterate via $x^{k+1}= x^k+\xi h^k$.
\EndWhile
\end{algorithmic}
\end{algorithm}

This two-level additive Schwarz operator (even with $\xi = 1$) can also be used as a preconditioner for the conjugate gradient method, which results in a substantial improvement in convergence. In the symmetric case, an error reduction factor of $2(1-2/(1+\sqrt \kappa))$ per iteration step is then obtained in contrast to the reduction factor of $1-2/(1+\kappa)$ from (\ref{redrate}). 

Note here again that, for general diagonal matrices $D_i$, the associated preconditioners are no longer symmetric, while they are in the case $D_i = \omega_i I$. This can cause both theoretical and practical problems for the conjugate gradient iteration. Then, instead of the conventional conjugate gradient method, we could resort to the {\em flexible} conjugate gradient method, which provably works also in the non-symmetric case, see  \cite{BDK15} and the references cited therein. 
This issue, however, can be completely avoided in our setting with the right choice of the overlap parameter $\gamma$.

\subsection{Outlook to fault-tolerance}
\label{subsection:outlook_fault_tolerance}
Now, let us shortly look ahead to our overall goal of developing a fault-tolerant domain decomposition solver for the combination method. In principle, we may of course choose any value for the overlap parameter $\gamma > 0$ and we may select quite freely the $D_i$ or $\omega_i$. However, if we want to achieve fault tolerance and thus need to employ a larger overlap of the subdomains, which are created as described above via our space-filling curve construction to allow for proper redundant storage and for data recovery, it is sensible to restrict ourselves to values of $\gamma$ that are integer multiples of $0.5$. In this case, every fine grid point is overlapped by exactly $2\gamma + 1$ subdomains, whereas, if $\gamma$ is not an integer multiple of $0.5$, the number of subdomains that overlap a particular point can indeed be different for different points of the same subdomain. Additionally, integer multiples of $0.5$ for $\gamma$ are the natural redundancy thresholds of our fault-tolerant space-filling curve algorithm. In particular for $\frac 1 2 n\leq\gamma < \frac 1 2 (n+1), n \in \mathbb{N}$, our fault-tolerant algorithm can recover from faults occurring for at most $n$ neighboring processors in the same iteration \cite{Griebel.Schweitzer.Troska:2020}. Thus, with these considerations, overlap parameter values of the form $\gamma=\frac 1 2 n, n \in  \mathbb{N},$ are the ones that are most relevant for proper redundant storage, for data recovery and thus for fault tolerance in practice.
Additionally, such a specific choice of $\gamma$ has a direct consequence on the resulting $\omega_i$ and $D_i$. According to \cite{Griebel.Schweitzer.Troska:2020} there holds the following Lemma which is merely a consequence of our particular construction of the overlaps and is valid for any space-filling curve.

\begin{lemma}\label{lemweig}
Let $d\ge1$ be arbitrary and let $\gamma = \frac{1}{2}n$, where $n\in\mathbb{N}, n\leq P - 1$. Then, with $c\coloneqq 2\gamma + 1$, there holds
\begin{equation*}
 D_i = \omega_i I = \frac{1}{c} I
\end{equation*}
for all $i=1,\ldots, P$ and any type of space filling curve employed.
\end{lemma}

Thus the general weightings $D_i$ and $\omega_i$ for the fine scale subdomains are the {\em same} and even boil down to a simple constant global scaling with the factor $1/c = 1/(2\gamma + 1) = 1/(n+1) $ uniformly for all subdomains $i=1,\ldots,P$, if the overlap parameter satisfies $\gamma = \frac{1}{2}n$, $n\in\mathbb{N}$. Thus, we obtain symmetry of the corresponding operators even for the $D_i$-choice. Note that this symmetry of the operator can not be obtained easily for the general $D_i$-weighting within the conventional geometric domain decomposition approach for $d>1$. Note also that the result of Lemma \ref{lemweig} holds analogously for more general geometries of $\Omega$ beyond the $d$-dimensional tensor product domain.
Note furthermore that, for the choice $\gamma = \frac{1}{2}n$, $n\in\mathbb{N}$, the weighted one-level operator (\ref{w1levprec}) becomes just a scaled version of the conventional one-level operator (\ref{AS2a}), i.e.
$$
C_{(1),D}^{-1}  =\frac 1 {2 \gamma+1}  C_{(1)}^{-1}.
$$
Consequently, we obtain
\begin{equation*}
C_{(2),D}^{-1} =  R_0^T A_0^{-1} R_0 + C_{(1),D}^{-1} = R_0^T A_0^{-1} R_0 + \frac 1 {2 \gamma+1}  C_{(1)}^{-1}.
\end{equation*}
and 
\begin{equation*}
C_{(2),D,bal}^{-1}= G^T C_{(1),D}^{-1} G + F =  \frac 1 {2 \gamma+1} G^T C_{(1)}^{-1} G + F,
\end{equation*}
which now just resemble a fine-level-rescaled variant of the conventional two-level operator and of its balanced version, respectively.


%
%
%
\section{Numerical Results}
\label{section:numerical_results}
We will consider the elliptic diffusion-type model problem 
$$
- \nabla \alpha ({\bf x}) \nabla u({\bf x}) = f({\bf x}) \quad \mbox{ in } \Omega=[0,1]^{d} 
$$
with right hand side $f({\bf x})$ and appropriate boundary conditions on $\partial \Omega$.
Since we are merely interested in the convergence behavior and the scaling properties of our approach and not so much in the solution itself, we resort to the simple
Laplace problem, i.e.~we set $\alpha= I$, $f=0$, and employ zero Dirichlet boundary conditions. Consequently, the solution is zero as well. For the discretization we employ finite differences on a grid with level parameter $l=(l_1,\ldots,l_d)$, which leads to $N$ interior grid points and thus $N$ degrees of freedom, compare (\ref{dof}), and which results in the associated matrix $A$.  Now any approximation $x^k$ during an iterative process directly gives the respective error in each iteration. We measure the error in the discrete energy norm associated to the matrix $A$ that stems from the finite difference discretization, i.e.~we track 
$$\lVert x^k\rVert_A := \sqrt{(x^k)^TAx^k}$$
for each iteration of the considered methods. Note here that we run the iterative algorithms for the symmetrically transformed linear system ${\hat A \hat x= \hat b}$ with ${\hat A=T^TAT}$, ${\hat b= T^Tb}$, ${\hat x=T^{-1}x}$ and ${T= diag(A)^{-1/2}}$,  whereas we measure the error in the untransformed representation, i.e.~for $x-x^k$. 
For the initial iterate $x^0$ we uniformly at random select the entries ${\tilde x_i^0, i=1\ldots,N,}$ of $\tilde x^0 $ from $[-1,1]$ and rescale them via $x_i^0 := \tilde x_i^0/\lVert\tilde x^0\rVert_A$ such that $\lVert x^0\rVert_A=1$ holds. To this end, we employed the routine \texttt{uniform\_real\_distribution} of the C++ STL (Standard Template Library).
We then run our different domain decomposition solvers until the relative error in the energy norm is reduced by at least a factor $10^{-8}$ and record the necessary number of iterations $K$, i.e.~$\lVert x^K\rVert_A \leq 10^{-8}$.\footnote{We confirmed via further experiments that all provided iterations showed a rather short preasymptotic phase so that the provided number of iterations are very much related to the asymptotic convergence rate of the respective iterative solver.} 

In the following, we present the results of our numerical experiments\footnote{All 
calculations have been performed on the parallel system 
{\em Drachenfels} of the Fraunhofer Institute for Algorithms and Scientific Computing (SCAI).
It provides, among others, $1.824$ Intel Sandy Bridge cores on $114$ compute nodes, each one with $2$ Xeon E5-2660 processors ($8$ cores per CPU, disabled SMT) at $2.20$ GHz and $32$ GB RAM, i.e.~2GB RAM per core, and  $2.272$ Ivy Bridge (Xeon E5-2650 v2) cores on $142$ compute nodes, each one with $2$ Xeon E5-2650 processors ($8$ cores per CPU, disabled SMT) at $2.60$ GHz and $64$ GB RAM, i.e.~4GB RAM per core. Drachenfels is equipped with a Mellanox Technologies MT27500 Family [ConnectX-3] 4x (56Gbps) connection. 
} 
using the additive and balanced two-level preconditioners given by \eqref{addSW} and \eqref{balSW}, respectively, each in their unscaled, $\omega_i$-scaled \eqref{omegai} and $D_i$-scaled \eqref{Di} variant. These operators are used within a Richardson-type iteration as given by Algorithm \ref{Alg:Richardson} and the correspondingly preconditioned conjugate gradient method. For the Richardson-type iteration we employ the optimal damping parameter $\xi^*=2/(\lambda_{min}+\lambda_{max})$ to obtain the best possible solver performance for comparison, even though such a damping is typically infeasible as it requires to determine the two extremal eigenvalues of the respective preconditioned operator by some numerical method beforehand. Furthermore, we present first scalability results within the context of a combination technique simulation in dimensions $d=2,3,4,5,6$ utilizing up to one million processors. Recall that the overall method comprises two levels of parallelism, the first level being the independent subproblems of the combination technique and the second level the domain decomposition of each of these subproblems. Moreover, the number of independent subproblems in the combination technqiue grows very fast with increasing dimension, compare \ref{NSPr}, so that it is sufficient to consider scalability in the second level of parallelism, i.e. the number of subdomains in the solution of each subproblem, up to a rather moderate number of processes while attaining very large scalability for the overall method, see section \ref{section:application_in_ct}. 

\subsection{Convergence properties}
First, we consider the one-dimen\-sio\-nal case in our numerical experiments. Surely there is no need to employ a parallel iterative domain decomposition solver and a direct sequential Gaussian elimination would be sufficient. However, this is a good starting point to study the convergence and parallel scaling properties of the various algorithmic variants. 
Moreover it will turn out that the one-dimensional case is indeed the most difficult one for good convergence and scaling results.
This behavior stems from the relative ''distance'' of the fine scale to the coarse scale, ${\frac{H}{h}\approx(\frac{qP}{N})^{-1/d}}$, since from conventional domain decomposition theory we have (for geometric coarse grid problems) a condition number of the order $O(1 + H/h)$ (see (\ref{AAA}) with $\delta = ch$), which is maximal for $d = 1$.

\begin{figure}[tb]
\centering
	\includegraphics[width=0.82\textwidth,height=0.33\textheight]{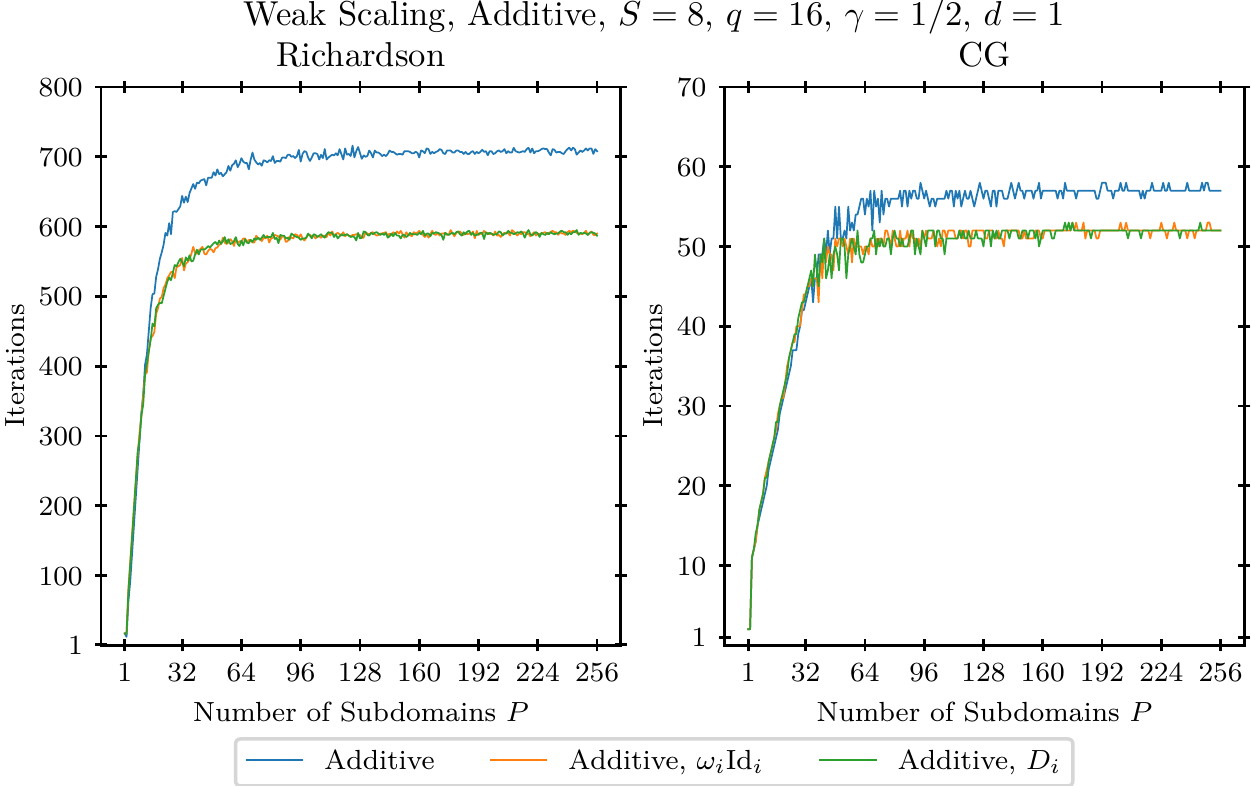}
	\includegraphics[width=0.82\textwidth,height=0.33\textheight]{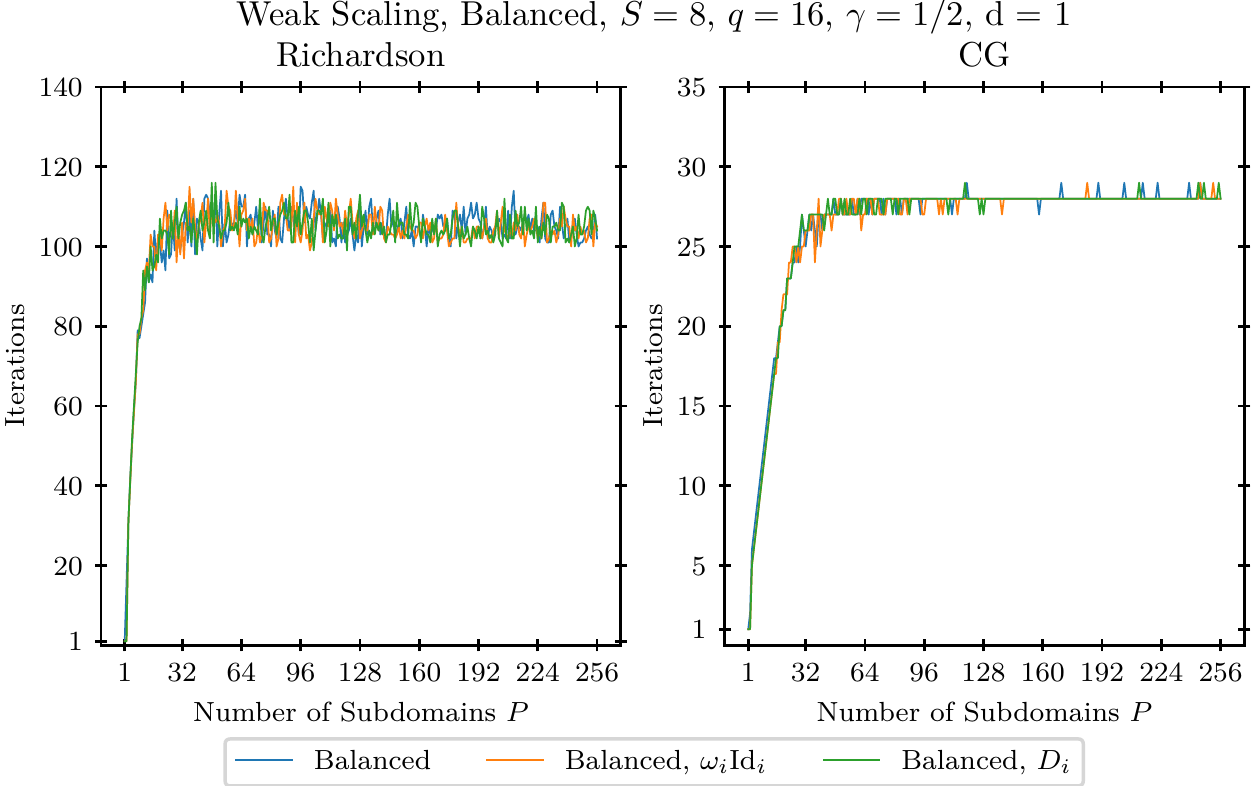}
	\caption{Weak scaling: Number of iterations versus number of subdomains for the linear/Richardson iteration (left) and the associated preconditioned conjugate gradient method (right), preconditioned by the three differently weighted versions of the {\em additive} (top) and the {\em balanced} (bottom) scheme, $d=1$, $q=16$, $\gamma=0.5$.}
\label{compare1_it}
\label{compare2_it}
\end{figure}

We set $N:= 2^S P$ where $S:=8$. Thus the size of each subproblem stays fixed at $2^8$ with growing $P$, whereas the overall number $N$ of unknowns grows linearly with $P$. Moreover we fix $q=16$ and $\gamma=0.5$. We compare the different methods for the three scalings $D_i=I$ (no scaling at all), $D_i=\omega_i\textrm{Id}_i$ with $\omega_i$ according to (\ref{omegai}), and $D_i$ according to (\ref{Di}), $i=1,\ldots,P$.
For our special choice of $\gamma$ we know from Lemma \ref{lemweig} that the weighting with $\omega_i$ and the weighting with $D_i$ are indeed the same and differ from the unweighted case by just the constant scaling $\frac 1 c I$ with $c=2 \gamma+1$. For the coarse scale problem, we always set $D_0=I$. 
Figure \ref{compare1_it} gives the number of iterations necessary to reduce the initial error by a factor of $10^{-8}$ 
for the Richardson-type approach (\ref{Alg:Richardson}) and the conjugate gradient method, both preconditioned via the additive scheme (\ref{addSW}) and the balanced variants according to (\ref{balSW}).
We obtain weak scaling behavior in all cases, i.e.~the necessary number of iterations stays constant for growing values of $P$. This constant depends on the respective splitting:
For the linear/Richardson iteration (left) we see that a scaling with $\omega_i$ (and equally with $D_i$) reduces this constant compared to the no-scaling case, albeit a large number of iterations is still needed. Moreover the results are substantially improved by the balanced variants: All three scalings now give the same results and the weak scaleup constant is reduced by a factor of approximately seven for the unweighted case and by a factor of approximately six for the other two variants.
For the associated preconditioned conjugate gradient methods (right) we observe a further reduction of the necessary number of iterations. 
This reflects roughly the $\kappa$-versus-$\sqrt \kappa$ effect of the conjugate gradient method in its convergence rate.
For the balanced version, we additionally see a substantial improvement of the scaleup constant compared to the additively preconditioned conjugate gradient method by a factor of nearly one half.
Note here that for values of $\gamma$ that are not integer multiples of 1/2, the $D_i$-scaling does not lead to a symmetric operator, which renders a sound and robust convergence theoretically questionable and it gave considerably worse iteration numbers with oscillating behavior for the corresponding preconditioned conjugate gradient method in further experiments.

We conclude that the Richardson iteration using the balanced preconditioner and the corresponding preconditioned conjugate gradient method are substantially faster compared to the additive version. We also see that, for our choice $\gamma=1/2$, balancing eliminates the difference of the unscaled and the $\omega_i$-scaled (and $D_i$-scaled) cases. Moreover the preconditioned conjugate gradient version is nearly quadratically faster and gives good weak scaling constants. 
Therefore, we will from now on focus on the optimally damped, {\em balanced} Schwarz/Richardson-type iteration 
as well as the correspondingly preconditioned conjugate gradient method. The type of damping we will choose, i.e.~none at all, $\omega_i$ according to (\ref{omegai}) or $D_i$ according to (\ref{Di}) is still to be determined. We refrain from employing the $D_i$-weighting for arbitrary choices of $\gamma$, since in general this results in a non-symmetric operator for which our theory is not valid anymore. 
It now remains to study the behavior of the unweighted and the $\omega_i$-weighted algorithms in more detail.

\begin{figure}[tb]
\centering
	\includegraphics[clip, trim=0 0.38in 0 0, width=0.82\textwidth]{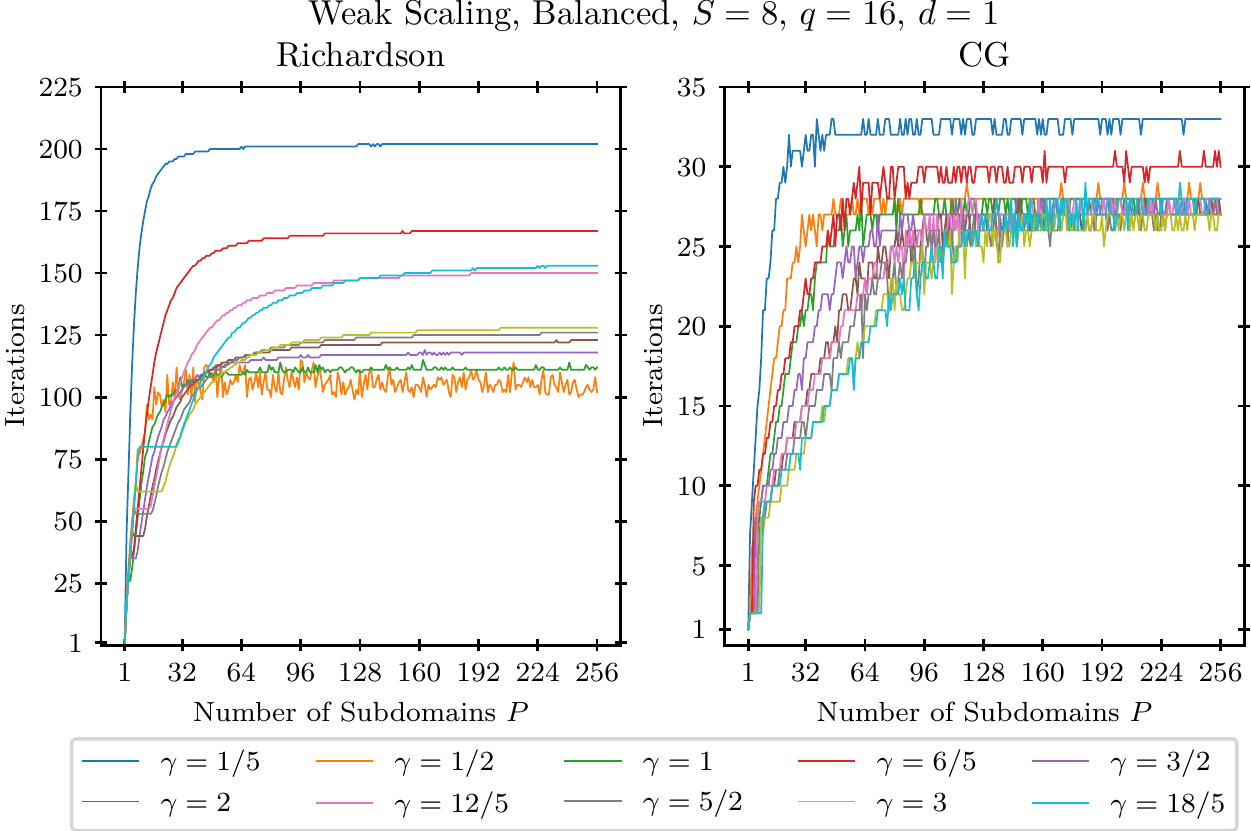}
	\includegraphics[width=0.82\textwidth]{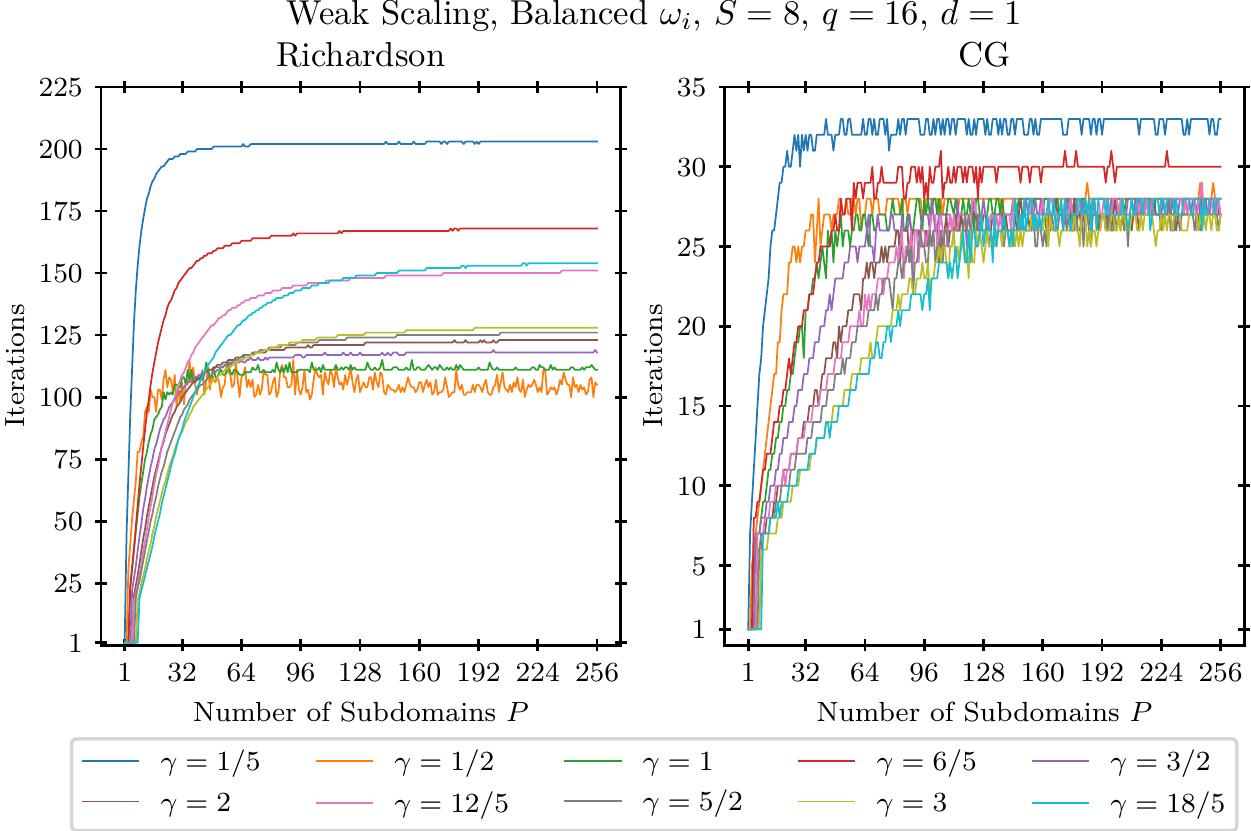}
	\caption{Weak scaling: Number of iterations versus number of subdomains for the unweighted balanced (top) and for the $\omega_i$-weighted balanced (bottom) linear/Richardson iteration (left) and the associated preconditioned conjugate gradient method (right) for different values of $\gamma$ with varying $P$, $d=1$, $q=16$, $S=8$.}
	\label{weak_vary_gamma_1d}
	\label{weak_vary_gamma_omegai_1d}
\end{figure}

So far, we kept the value of the overlap parameter $\gamma$ fixed. Now we vary $\gamma$ and consider its influence on the weak scaling behavior of our two algorithms in the balanced case. First we consider the one-dimensional situation, where we set $S=8$, $q=16$ and vary the number $P$ of subdomains. The resulting number of iterations for different values of $\gamma$ ranging from $1/5$ up to $5$ are shown in Figure \ref{weak_vary_gamma_1d}.
Comparing the $\omega_i$-weighted case to the unweighted case, there is not much visible difference at all.
We again clearly see weak scaling behavior.
The scaling constant now depends on the respective value of $\gamma$. 
In the Richardson case, it is interesting to observe that, in any case, a constant number of iterations is quickly reached for rising values of $P$. 
Moreover, for a small overlap value of $\gamma=1/5$, it is quite bad. The number of iterations is seen to be optimal for $\gamma=1/2$ and then deteriorates for larger values of $\gamma$. 
Note at this point that, starting with $\gamma=1/2$, we observe a slight deterioration of the convergence curves and of the weak scaling constant, where this deterioration is {\em monotone} in $n$ if we restrict ourselves to values of $\gamma$ that are integer multiples of $0.5$, i.e.~$\gamma=\frac 1 2 n, n \in \mathbb{N}_+$. The non-integer multiples give worse results.
In the conjugate gradient case, the scaling constant is reached increasingly later for rising values of $P$ (which is desirable). It is improved in a nearly monotonic way for rising values of $\gamma$. 
Furthermore the absolute value of necessary iterations is again much smaller than in the Richardson case.  

\begin{figure}[ht]
\centering
	\includegraphics[clip, trim=0 0.38in 0 0, width=0.82\textwidth]{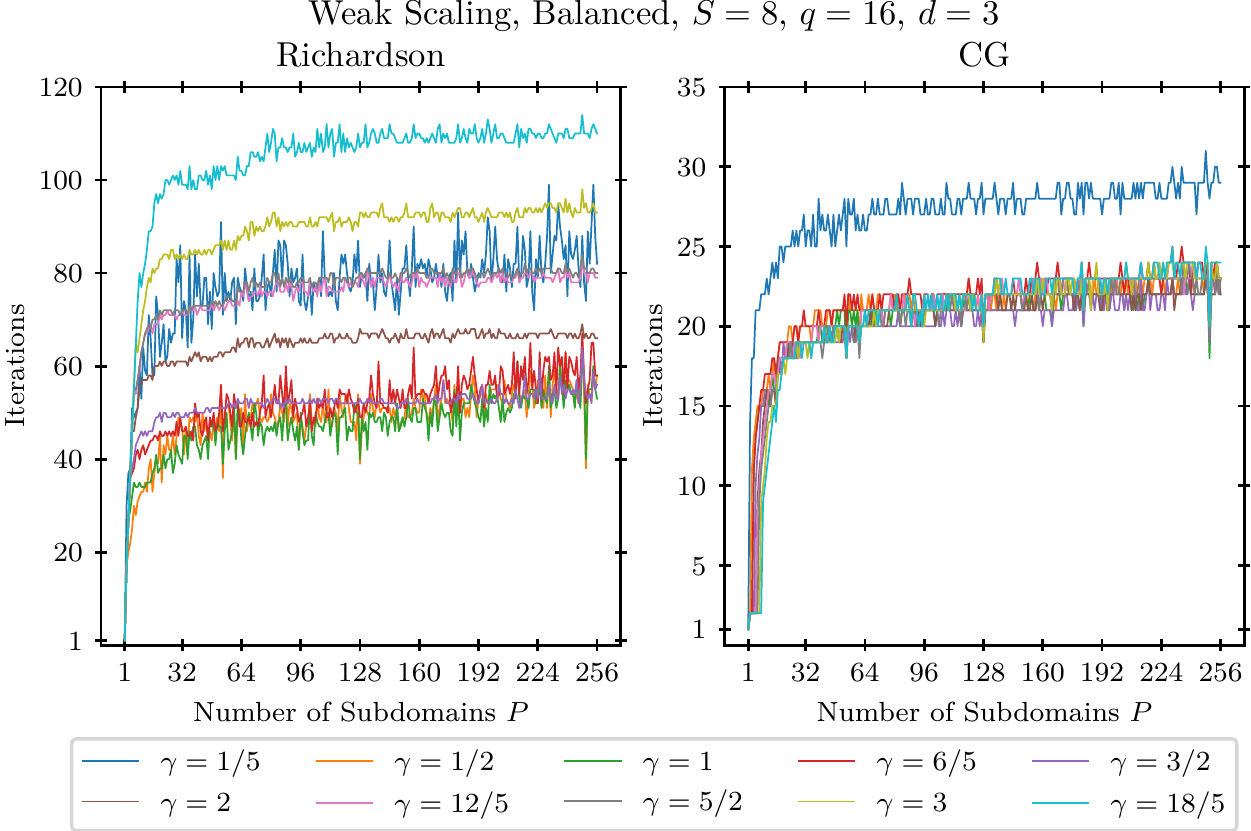}
	\includegraphics[width=0.82\textwidth]{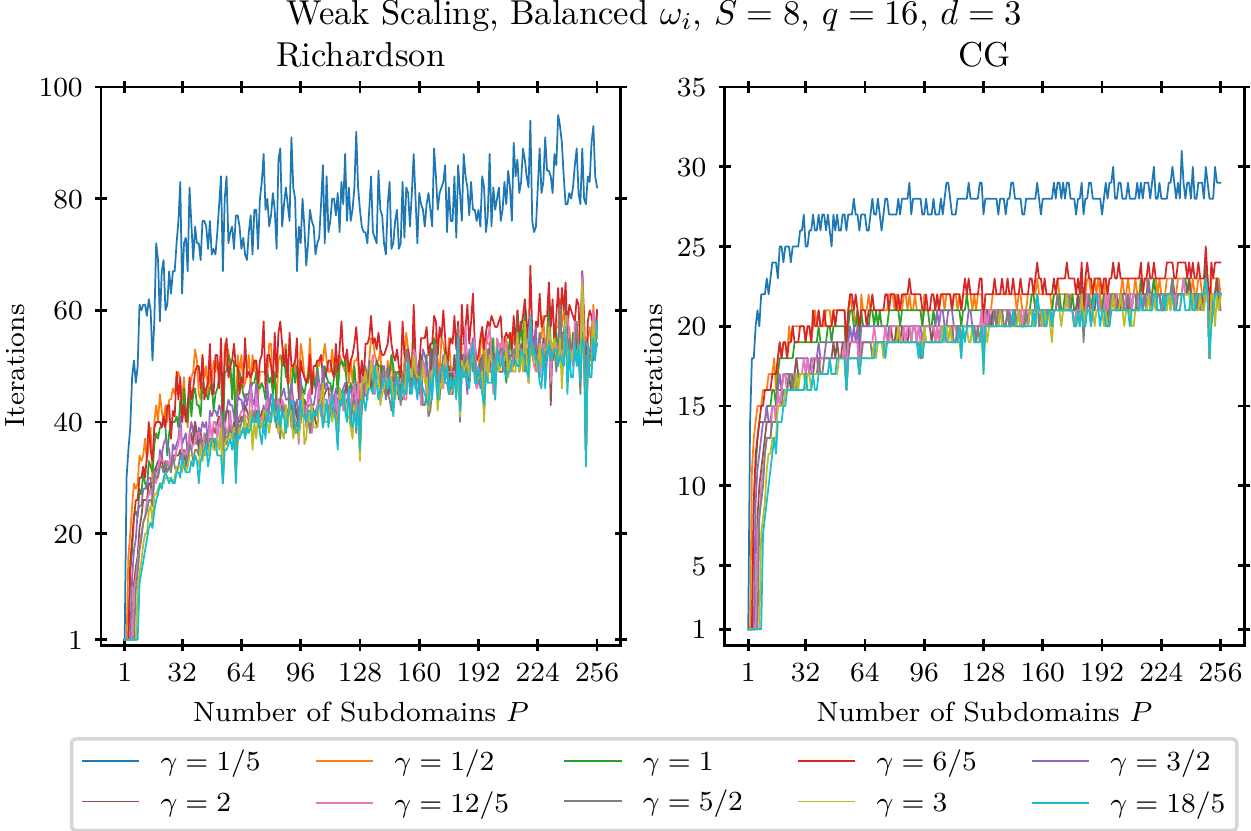}
	\caption{Weak scaling: Number of iterations versus number of subdomains for the unweighted balanced (top) and for the $\omega_i$-weighted balanced (bottom) linear/Richardson iteration (left) and the associated preconditioned conjugate gradient method (right) for different values of $\gamma$ with varying $P$, $d=3$, $q=16$, $S=8$.}
	\label{weak_vary_gamma_3d}
	\label{weak_vary_gamma_omegai_3d}
\end{figure}

Next, we consider the three-dimensional case and again vary the value of $\gamma$. Then, in contrast to the one-dimensional situation, the Hilbert curve structure comes into play. The resulting iteration numbers for different values of $\gamma$ are given in Figure \ref{weak_vary_gamma_3d}.

For the Richardson iteration we again see an analogous behavior for the weak scaling constant but with a much lower number of iterations compared to the one-dimensional case. Furthermore we observe that the $\omega_i$-scaling improves the convergence: Except for $\gamma=1/5$, all observed numbers of iterations are now close together for both the values of $\gamma$ that are integer multiples of $1/2$ and the other values, and their respective number of iterations are in any case reduced to around $50$. 
This shows that $\omega_i$-weighting is able to deal with the decomposition based on the Hilbert curve, which appears for $d>1$, in a proper way.
In the unweighted conjugate gradient case, the curves are approximately the same for all values of $\gamma$, except for $\gamma=1/5$ which is too small again. 
They are successively improved for rising values of $\gamma$ in the $\omega_i$-weighted case, as is expected intuitively. Again, the conjugate gradient method is much faster than the Richardson scheme. Similar observations could be made for other dimensions.

Altogether, $\omega_i$-weighting stabilizes the iteration numbers against variations of $\gamma$ and improves the convergence behavior. In light of the larger costs involved for higher values of $\gamma$, a good choice for $\gamma$ is given by the value $1/2$ in both the Richardson and the conjugate gradient case. 
However, larger values of $\gamma$ may be needed to attain fault tolerance due to redundancy later on.
The $\omega_i$-weighted case then indeed results in slightly better convergence results for larger values of $\gamma$ and shows no deterioration, which the unweighted case does.    
From now on we therefore will only consider the $\omega_i$-weighting for both schemes. Moreover we  will restrict ourselves in this work to the value $\gamma=1/2$. In view of fault-tolerance this value might need to be increased.

\begin{figure}[tb]
\centering
	\includegraphics[clip, trim=0 0.35in 0 0, width=0.82\textwidth]{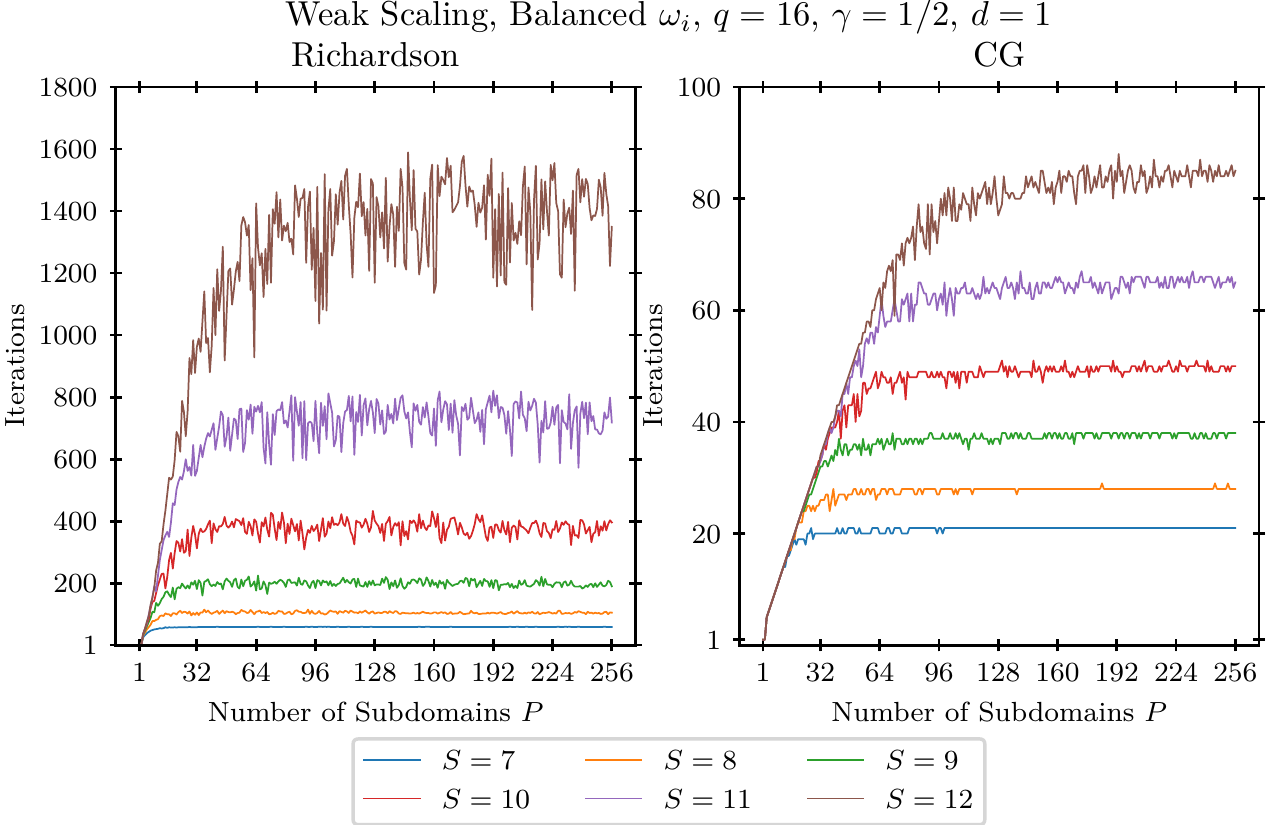}
	\includegraphics[width=0.82\textwidth]{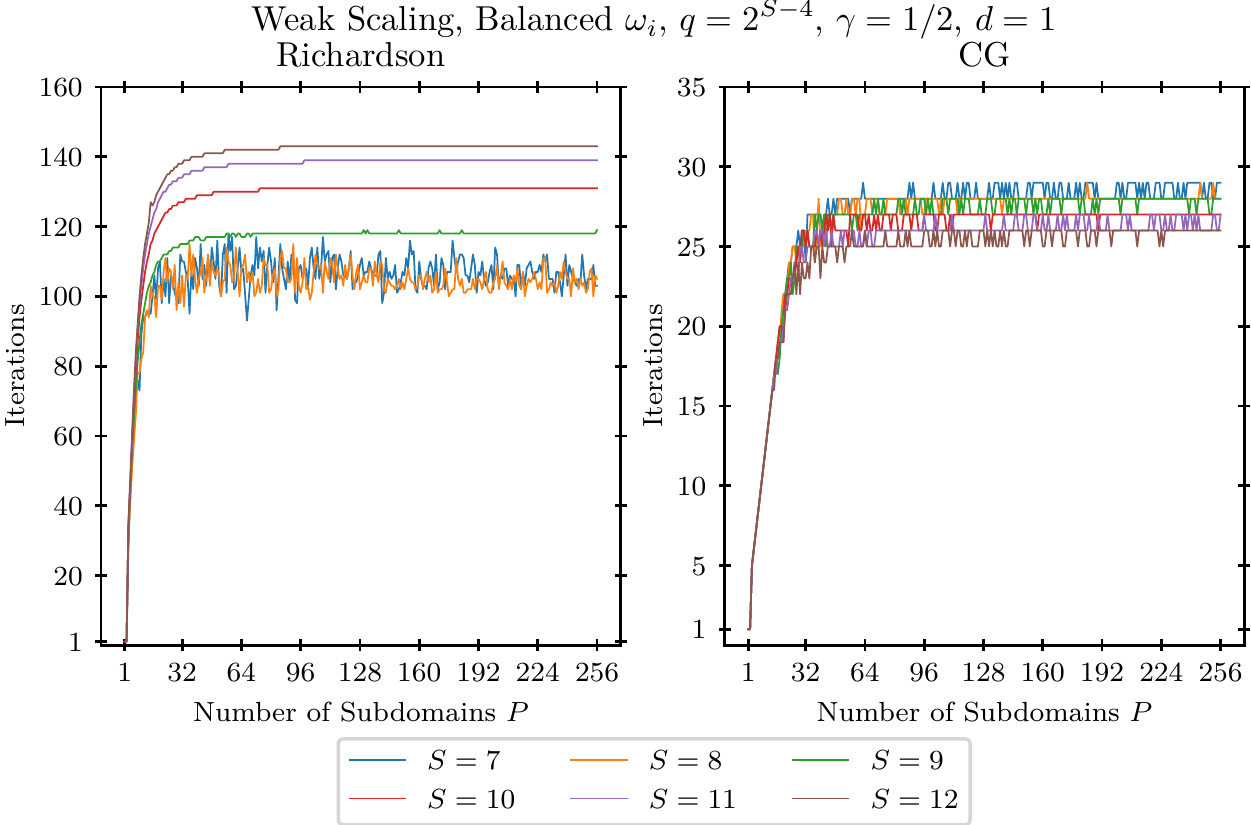}
	\caption{Weak scaling: Number of iterations versus number of subdomains for the $\omega_i$-weighted balanced linear/Richardson iteration (left) and the associated preconditioned conjugate gradient method (right) for different values of $S$, $d=1$ $\gamma=0.5$, $q=16$, (top) and $q=2^{S-4}$ (bottom).}
\label{fixq}
\label{varq}
\end{figure}

So far, we kept the weak scaling parameter $S$, i.e.~the size $2^S$ of each subproblem, fixed and only varied the number $P$ of subdomains. But what happens if we also vary the subproblem size? For fixed $q=16$, the results are shown in Figure \ref{fixq}.

We again clearly see weak scaling behavior.
However the weak scaling constant now depends on the subproblem size, i.e.~it grows with rising $S$. This holds for both the Richardson-type iteration and the conjugate gradient approach. This behavior stems from the fixed value of $q$ and thus the fixed size of the coarse scale problem for fixed $P$. In this case the difference between fine scale and coarse scale increases with growing $S$ and, consequently, the additive coarse scale correction is weakened compared to the fine scale.

We now let the coarse scale parameter $q$ be dependent on $S$. To be precise, we set $q=2^{S-4}$, which gives a coarse problem of size $P \cdot 2^{S-4}$, where we vary the parameter $S$. This way, we double $q$ while doubling the subdomain size $2^S$, i.e.~$N/P$.
The obtained results are shown in Figure \ref{varq} (bottom). 

In all cases, we obtain substantially improved results compared to the fixed choice ${q=16}$ from Figure \ref{fixq} (top). This was to be expected since now the coarse scale correction is improved for rising $S$ due to the constant relative ``distance'' between the coarse and fine scale. We again observe an asymptotically constant number of iterations for growing values of $P$ and we obtain a weak scaling constant which is, compared to the fixed choice ${q=16}$, now only slightly growing with $S$ for the Richardson iteration and seems to approach a limit of about 145 for rising values of $S$. Moreover it is now completely independent of $S$ for the conjugate gradient approach, for which we need at most 29 iterations in all cases. This shows that $q$ should scale proportional to $N/P$. In further experiments we quadrupled $q$ while doubling $N/P$ and found that the weak scaling constant then even shrank with growing $S$. 

Note at this point that the size of the coarse scale problem is now $2^{-4} \cdot P \cdot 2^S$, while the size of each subdomain problem is $2^S$ (whereas the size of the overall problem is $N=P\cdot 2^S$).
Thus the cost of solving the coarse scale problem tends to dominate the overall cost with rising $P$, which is the price to pay for a more uniform convergence behavior. This calls for further parallelization of the coarse scale problem itself via our $P$ processors to remedy this issue. Thus, in contrast to the present implementation via the Bank-Holst paradigm where we redundantly keep the coarse grid problem on each processor (besides the associated subdomain problem), we should partition the coarse scale matrix to the $P$ processors and therefore solve the coarse scale problem in a parallel way.
This however will be future work.
\begin{figure}[tb]
	\centering
		\includegraphics[width=0.82\textwidth]{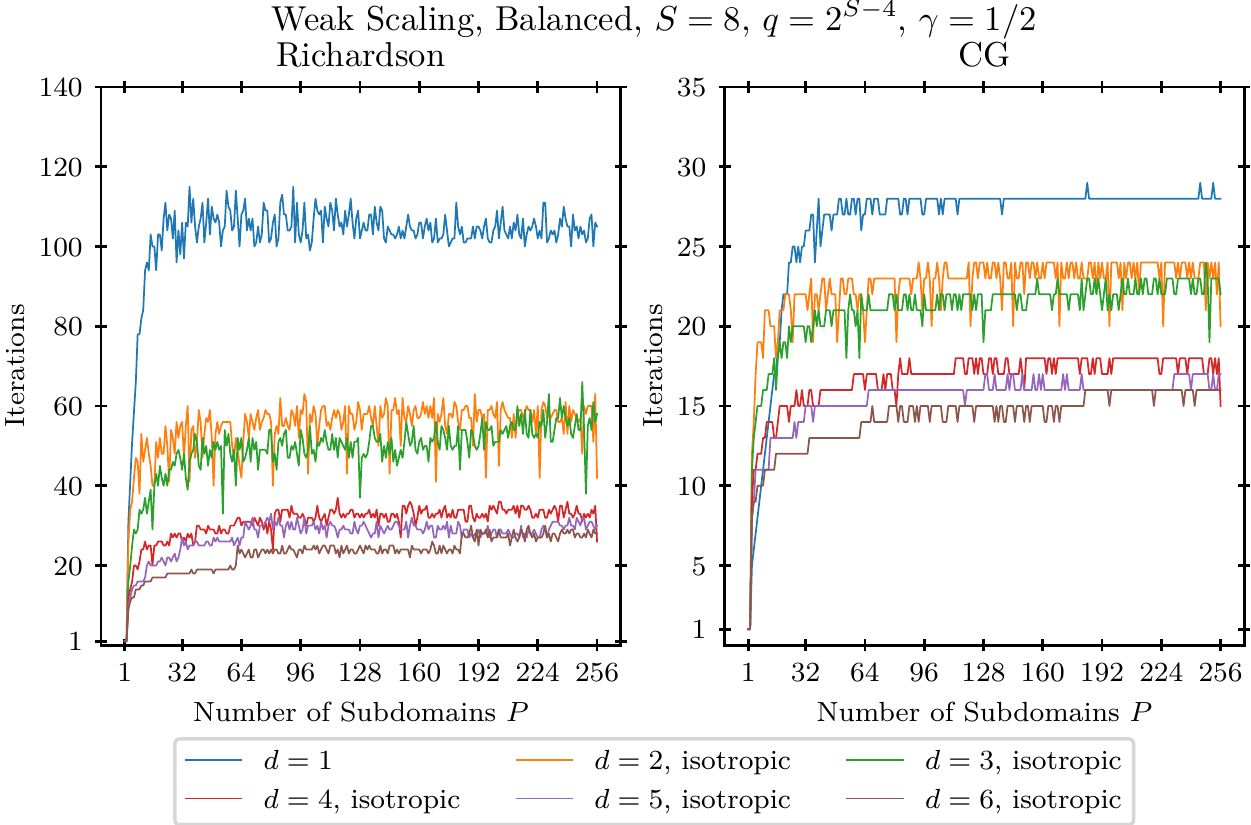}
		\caption{Weak scaling: Number of iterations versus number of subdomains for the $\omega_i$-weighted balanced linear/Richardson iteration (left) and the associated  preconditioned conjugate gradient method (right) for different dimensions $d$, $q=2^{S-4}$, $\gamma=0.5$.}
	\label{dimvarq}
	\end{figure}
	
\subsection{Parallel performance properties: Scaling behavior}
Next, we consider the weak scaling behavior for varying dimensions $d=1,2,3,4,5,6$.
For the discretization we stick to the isotropic situation, i.e.~we set
$${l} =(\lfloor (S+\log_2(P))/d\rfloor, \ldots, \lfloor (S+\log_2(P))/d \rfloor).$$ 
Thus the overall number of degrees of freedom is independent of $d$, since 
$$N \approx \prod_{j=1}^d 2^{(S+\log_2(P))/d} =2^{S+\log_2(P)}=P \cdot 2^S,$$ 
and the size of each subdomain is again approximately $2^S$ for all values of $P$. Furthermore we choose $\gamma=0.5$, $S=8$ and $q=2^{S-4}$ and consider only the $\omega_i$-weighted balanced methods. The resulting weak scaling behavior is shown in Figure \ref{dimvarq}.

We always obtain weak scaling behavior, where now the constant depends on $d$. But it improves for growing values of $d$ and we see that the one-dimensional case is indeed the most difficult one. This is due to the relative ''distance'' of
the fine scale to the coarse scale in the two-level domain decomposition method, which is largest for $d=1$  and decreases for larger values of $d$. Furthermore a stable limit of 26 and 16 iterations, respectively, is reached  for $d\geq 6$.
Such a behavior could be observed not only for the isotropic discretizations in $d$ dimensions but also for all the various anisotropic discretizations, which arise from (\ref{eq:combi}) in the sparse grid combination method. This becomes clear when we consider the simple case of the anisotropic discretization $l=(\lfloor S+\log_2(P)\rfloor,1,\ldots, 1)$ in $d$ dimensions: With our homogeneous Dirichlet boundary conditions we obtain the same fine scale finite difference discretization as for the case $d=1$ with differential operator $-\partial^2/\partial^2_{x_1} + 8 (d-1) \cdot I_{x_1}$. We then have an additional reaction term of size $8 (d-1),$ which merely improves the condition number compared to the purely one-dimensional Laplacian. Consequently, the one-dimensional convergence results impose an upper limit to the number of iterations needed for all the subproblems arising in the combination method.
	
\begin{figure}[tb]
	\centering
		\includegraphics[clip, trim=0 0.2in 0 0, width=0.82\textwidth]{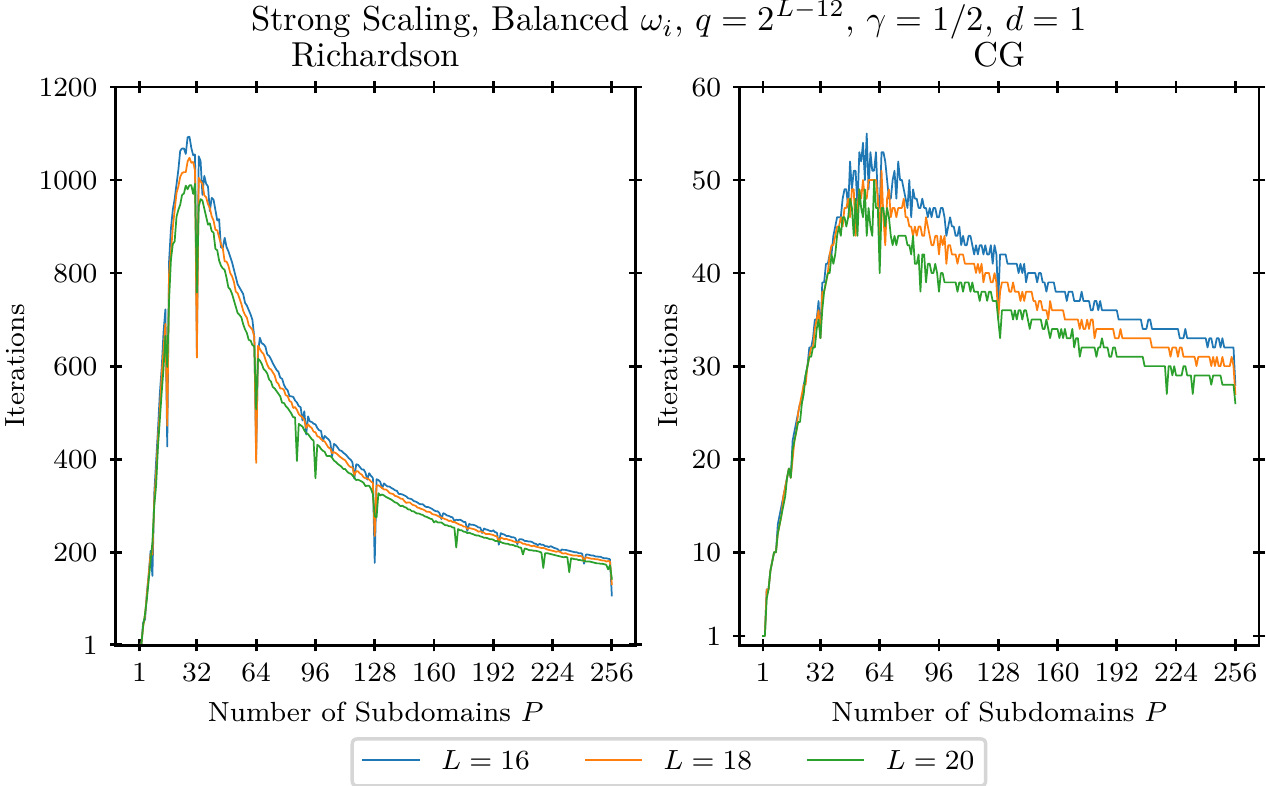}
		\includegraphics[width=0.82\textwidth]{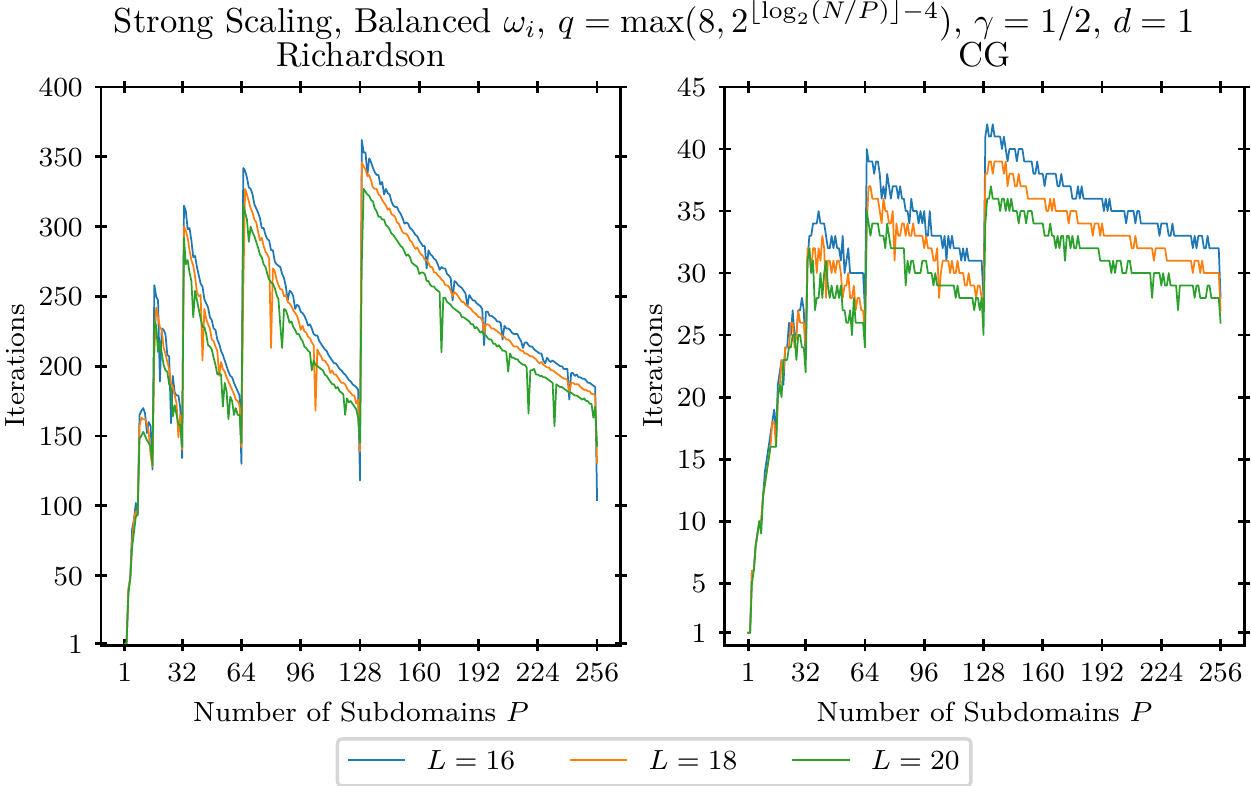}
		\caption{Strong scaling: Number of iterations versus number of subdomains for the $\omega_i$-weighted balanced linear/Richardson iteration (left) and the associated preconditioned conjugate gradient method (right), $d=1$, ${L=16,18,20}$, $\gamma=0.5$, $q=2^{L-12}$ (top) and ${q= 2^{\lfloor \log_2 (N/P) \rfloor -4}}$ (bottom).}
	\label{strongvar}
	\label{strongfix}
	\end{figure}

Now let us shortly consider the strong scaling situation as well. There we have $N=2^L$, where the size of each subdomain is $2^L /P$, i.e.~it decreases with growing values of $P$. Moreover values of $P$ larger than $2^L$ are not feasible. We consider again the one-dimensional situation, set $q=2^{L-12}$, $\gamma=0.5$, and vary the number $P$ of subdomains. The resulting strong scaling behavior is shown for $L=16,18,20$ in Figure \ref{strongfix}.

We see that the necessary number of iterations first grows with rising values of $P$ and then, after its peak, steadily declines, as is expected. This is due to the fact that now the size of each subproblem shrinks with rising $P$, whereas the size $P \cdot 2^{L-12}$ of the coarse scale problem grows linearly with $P$. There is not much of a difference for the three curves $L=16$, $L=18$ and $L=20$ since  for each $P$ the coarse scale problem has the same relative distance to the fine grid for all values of $L$, i.e.~$ 
2^{16}/(2^{16-12} \cdot P) = 2^{18}/(2^{18-12} \cdot P)=2^{20}/(2^{20-12} \cdot P)$. Note here that the downward peaks of the curves in Figure \ref{strongfix} (top) correspond to the situation
where all subdomains are perfectly balanced, i.e.~where $N \mod P=0$.  
Furthermore the parameter $q$ was chosen such that the case $P=256$ for $L=16,18,20$ in Figure \ref{strongfix} (top) results in exactly the same situation as the case $P=256$ for $S=8,10,12$ in Figure \ref{varq}, respectively.

We finally consider strong scaling for the situation where $q= 2^{\lfloor \log_2(N/P)\rfloor-4}$, i.e.~where the choice of $q$ is $P$-dependent and is chosen such that the relative distance between the coarse and the fine grid is roughly four levels for all $P$, as was the case in the weak scaling experiments in Figure \ref{varq}. The results are given in Figure \ref{strongvar} (bottom). Due to the involved round off there are now jumps at the points where $\log_2(P)$ is an integer. 
However we observe the expected steady decline of the necessary number of iterations in between these jump points. Note that, for this choice of $q$, the necessary numbers of iterations for the end points of these intervals are approximately the same.


\subsection{Application within the combination technique}
\label{section:application_in_ct}
In the last experiment we demonstrate the use of our dimension-oblivious domain decomposition method within the combination technique. To this end, we consider the Poisson model problem
\begin{equation}
	\label{eqn:poisson:dim}
	- \Delta u({\bf x}) = f({\bf x}) \quad \mbox{ in } \Omega=[0,1]^{d} \quad \mbox{ for } d=2,\ldots, 6
\end{equation}
where we choose $f$ and the employed Dirichlet boundary conditions such that the exact solution is given by
$$
u({\bf x}) = \lVert{\bf x}\rVert_2\prod_{i=1}^d\sin(\pi x_i).
$$

We solve this problem on different discretization levels for the full grid as well as the combination technique using the conjugate gradient method preconditioned by our $\omega_i$-weighted balanced domain decomposition method. The stopping criterion of the solvers was chosen as the relative residual norm of $10^{-8}$, i.e.\ $\frac{\lVert r_k \rVert}{\lVert r_0 \rVert} < 10^{-8}$.

\begin{figure}[tb]
	\centering
	\includegraphics[width=0.82\textwidth]{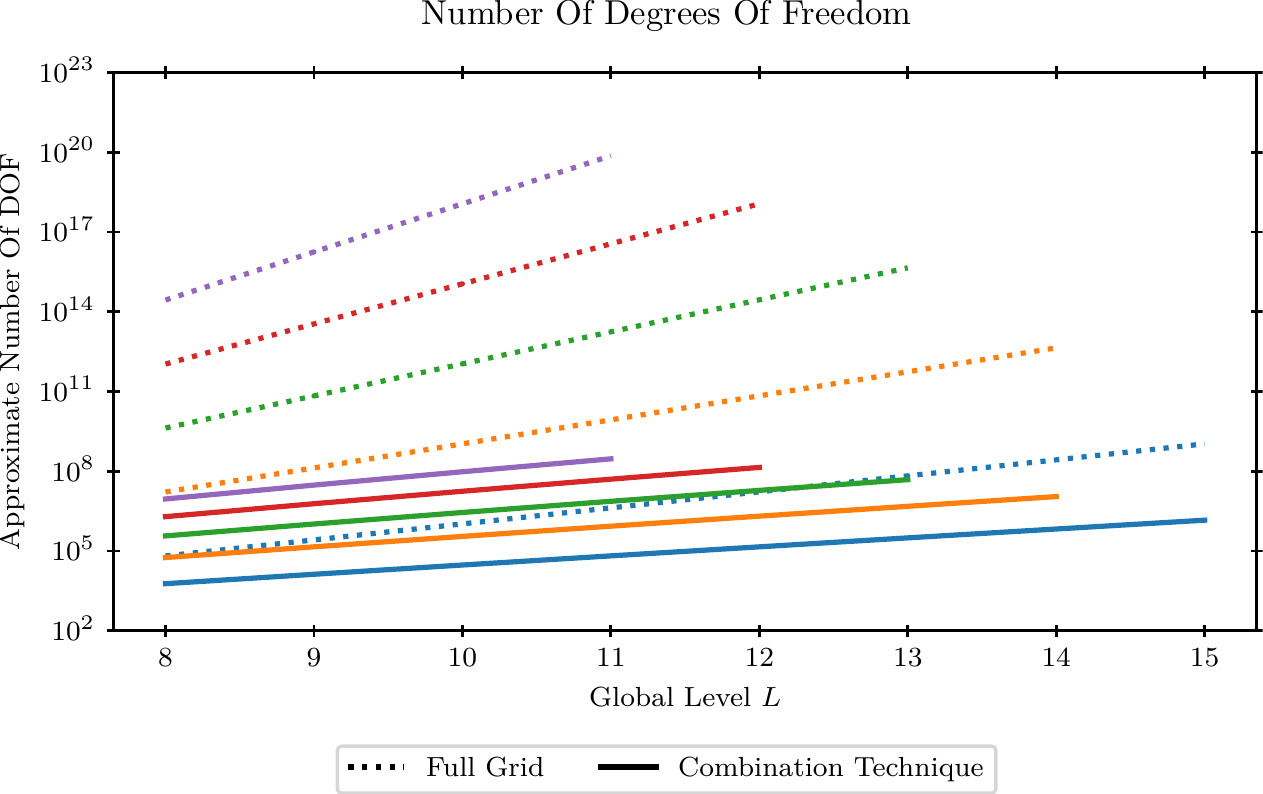}
	\includegraphics[width=0.82\textwidth]{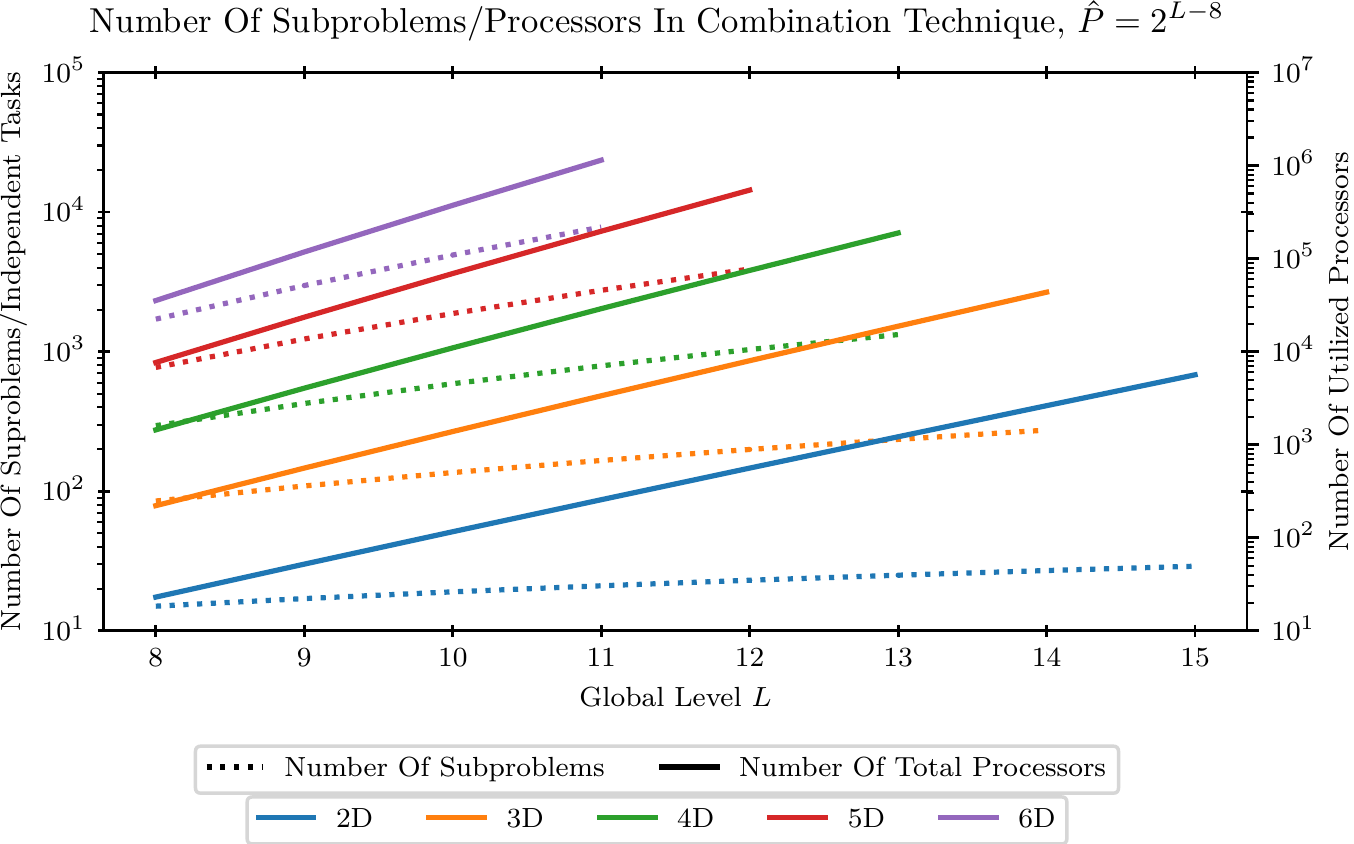}
	\caption{Approximate number of degrees of freedom for the combination technique and full grid (top) and number of subproblems and utilized processors in the combination technique, $\hat P = 2^{L-8}$ (bottom).}
\label{ct_results_dof_processor}
\end{figure}

Figure \ref{ct_results_dof_processor} (top) shows the approximate total number of degrees of freedom for both discretizations, where for the discretization parameter $L$ the full grid contains $2^L-1$ interior nodes in each coordinate direction and the combination technique is constructed as described in section \ref{section:motivation}. Clearly, due to the curse of dimensionality, the full grid discretization quickly becomes infeasible since the number of degrees of freedom scales as $O(2^{Ld})$. On the other hand, the total number of degrees of freedom for the sparse grid approach is approximately $O(2^L(\log(2^L))^{d-1})$ only. Moreover, the problem \eqref{eqn:poisson:dim} is never discretized explicitly employing the total number of degrees of freedom due to the combination technique. The problem is rather discretized on a large number of independent subproblems on significantly smaller and in general anisotropic classical grids, compare Figure \ref{fig:combi_scheme}.
The number of these combination technique subproblems for different discretizations can be found in Figure \ref{ct_results_dof_processor} (bottom). Each subproblem is fully independent of all other subproblems which results in a first level of parallelism. The second level of parallelism is achieved by solving each subproblem using our domain decomposition approach. Recalling the choice \eqref{Pchoice} of the number of subdomains/processors per subproblem via the speedup factor $\hat P$, the choice $\hat P = 2^{L-8}$ results in every subdomain of each subproblem containing approximately $2^8$ degrees of freedom, independent of the dimension $d$ and layer $i$ in the combination technique. The resulting equal amount of work is comparable to the weak scaling experiments of the previous section with $S=8$. Additionally, the finest level in dimension $d$ was selected such that it results in at most $256$ processors used for each subproblem. Hence, by definition \eqref{Pchoice}, we therefore require $2^8=256\ge P=\hat P \cdot 2^{d-1-i}=2^{L-8+d-1-i}$ for $i=0,\ldots,d-1$ and all $L$. Therefore the largest possible discretization level in dimension $d$ is $L_{\text{max},d}=17-d$. Furthermore, the coarse grid parameter $q$ for both the full grid and combination technique subproblems was selected such that the resulting coarse grid is four levels coarser then the fine grid, again in accordance to the results of the previous section. One can clearly see that with this moderate parallelism employed for the treatment of each subproblem the combination technique easily utilizes more than one million processors for the efficient solution of the overall global problem.  
	
\begin{figure}[tb]
	\centering
	\includegraphics[width=0.82\textwidth]{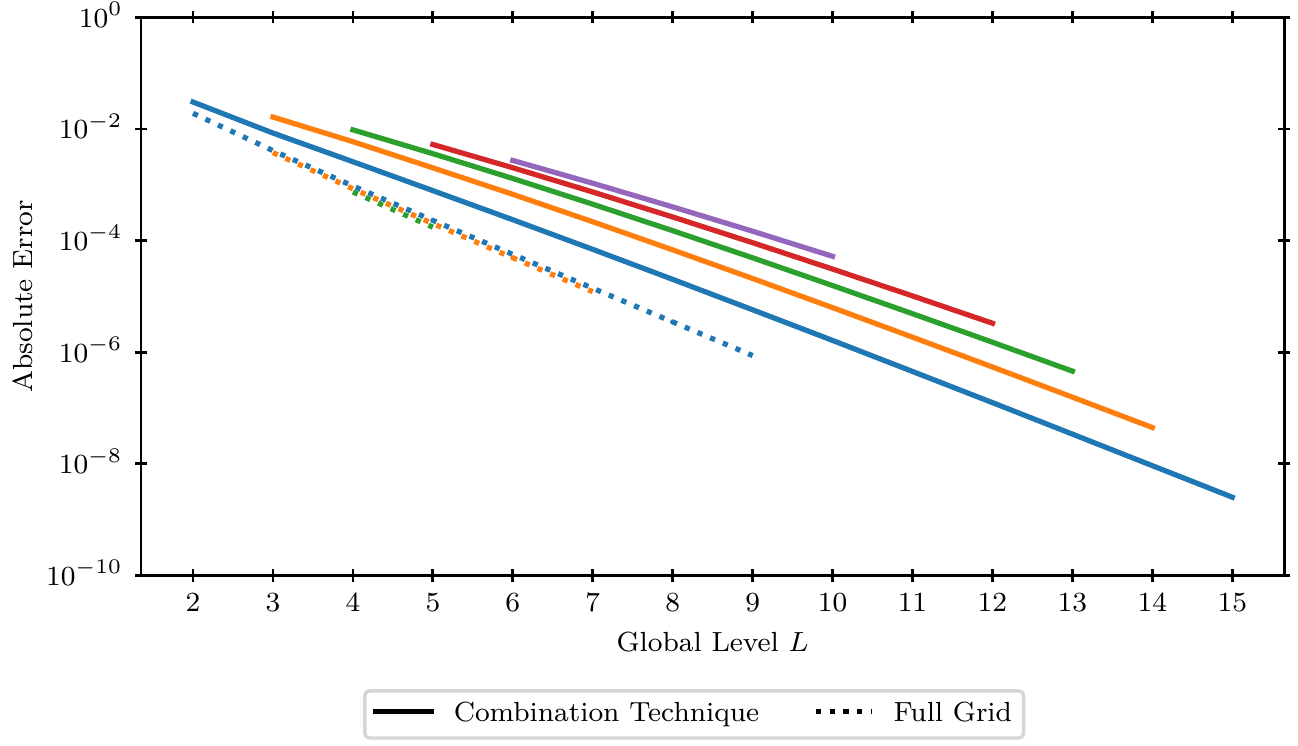}
	\includegraphics[width=0.82\textwidth]{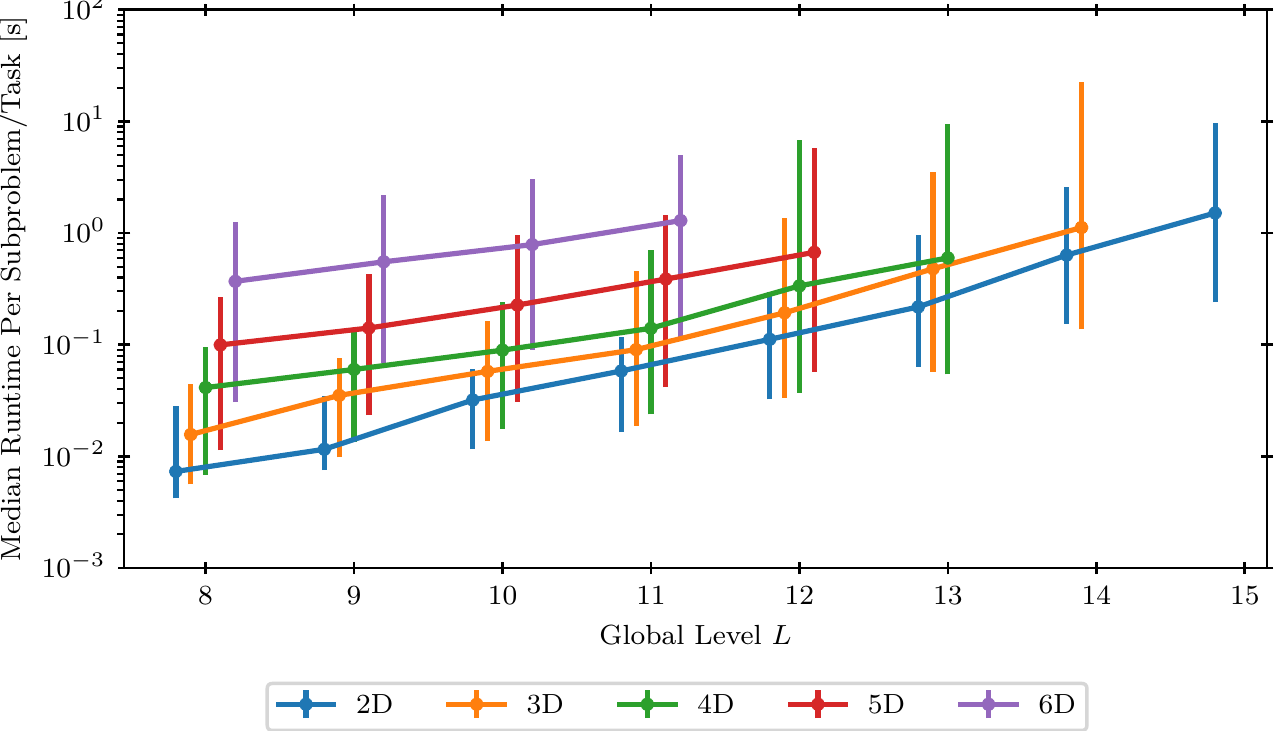}
	\caption{Absolute error of the solution for the full grid and combination technique (top) and runtimes for the combination technique subproblems (bottom) for different dimensions $d$ and discretization levels $L$.}
\label{ct_results_error_times}
\end{figure}

Figure \ref{ct_results_error_times} (top) shows the absolute error of both approaches. For the full grid we observe the theoretic quadratic convergence independent of the dimension. On the other hand, for the combination technique the theoretic convergence rate is given by $O(N^{-2}\log(N)^{d-1})$ with $N=2^L-1$. The dimension-dependent logarithmic term can clearly be observed. Note that only a small number of data points are available for the full grid compared to the sparse grid approach, which is again due to the exponential growth of the number of degrees of freedom and the consequent infeasibility of the problem. For the same reason, the error of the solution of the combination technique was only computed for a uniformly at random selected subset of the full grid nodes for combinations of the parameters $d$ and $L$ where no full grid data is available. Additionally, for $d=6$ and $L=11$ the error could not be computed since the full grid space-filling curve indices are not storable in the built-in unsigned integers, which are at most $64$-bit. However this does not affect the actual computation in the combination technique, which employs significantly smaller anisotropic grids, it only affects the post-processing of the error, which has been performed on the full grid for convenience. In fact, problems on much finer discretizations in higher dimensions than the ones considered can be solved using the combination technique without running into these issues.

The runtimes for each subproblem can be seen in Figure \ref{ct_results_error_times} (bottom). The plot shows the minimal, maximal as well as the median runtime of the subproblems. We can see that our choice of the parallelization parameter $\hat P$ yields scalability of the overall combination technique. Due to the robustness of the domain decomposition solver with respect to the dimension as observed in the previous section, we obtain an algorithm whose runtime depends on the number of degrees of freedom per processor only. This can clearly be observed for the finest discretizations in each dimension, where by construction each processor stores roughly $2^8$ degrees of freedom which yields very similar runtimes for the subproblems. The same can be observed for all combinations of $d$ and $L$ which result in the same number of processors per subproblem on the first layer, $P=2^{L+d-1}$, i.e.\ where $L+d$ yields the same value. Note again that the problem on the finest level for $d=6$ employs over a million processors, compared to approximately ten thousand processors for the finest level for $d=2$, while mainting similar runtimes for each subproblem, demonstrating the weak scalability of the overall algorithm independent of the problem dimension.


\section{Concluding Remarks}
\label{section:conclusion}
We presented an algebraic dimension-oblivous linear solver based on a space-filling curve domain decomposition approach with large overlap. The proposed solver moreover allows for the efficient utilization of arbitrary processor numbers while maintaining optimal convergence behavior. 

The presented domain decomposition linear solver is an essential ingredient in the design of an extremely scalable sparse grid combination technique solver for high-dimensional PDE problems of elliptic and parabolic type. Besides scalability results of the proposed solver we gave first results for the treatment of a simple diffusion problem in dimension $d=2,3,4,5,6$ via the combination technique and the proposed domain decomposition solver utilizing up to one million proceessors. The extension to fault-tolerance is currently under development and will be discussed in detail in a forthcoming paper, see however \cite{Griebel.Schweitzer.Troska:2020} for initial results.


\subsection*{Acknowledgment}
Michael Griebel was supported by a research grant 
of the {\em Bayerische Akademie der Wissenschaften} for the project {\em Stochastic subspace correction as a fault-tolerant sparse grid combination method on massive parallel compute systems for high-dimensional parametric diffusion problems}. The support of the Leibniz-Rechenzentrum in Garching, the support of Prof. Dr. Dieter Kranzlm\"uller, LRZ Garching, LMU M\"unchen, Institut f\"ur Informatik, and the support of Prof. Dr. Hans-Joachim Bungartz, TU M\"unchen, Institut f\"ur Informatik,
is greatly acknowledged. Michael Griebel thanks the Institut f\"ur Informatik of the LMU M\"unchen and the Leibniz-Rechenzentrum of the Bayerische Akademie der Wissenschaften for their hospitality.


\end{document}